\documentclass[12pt]{article}

\title{Non-locally modular regular types in classifiable theories}
\author{Elisabeth Bouscaren\thanks{Partially supported by ANR AAPG2019 GeoMod and NSF grant DMS-1855789.}\and Bradd Hart\thanks{Partially supported by NSERC.} \and Ehud Hrushovski \and
Michael C.\
Laskowski\thanks{Partially supported
by NSF grants DMS-1855789 and 2154101. }
}

\usepackage{amssymb, latexsym,amsmath}
\newbox\smilebox
\newbox\anchorbox
\newbox\noanchorbox
\newbox\tempbox

\setbox\smilebox=\hbox{$\smile$}

\def\anchor{\hbox{\vtop{
           \hbox to \wd\smilebox{\hfil\vrule width.4pt height7pt depth1pt\hfil}
           \vskip  -11.5truept
           \hbox to \wd\smilebox{\hfil$\smile$\hfil}}}}
\setbox\anchorbox=\anchor
\def\noanchor{\hbox{\vtop{
           \hbox to \wd\anchorbox{\hfil\anchor\hfil}
           \vskip -14truept
           \hbox to \wd\anchorbox{\hfil/\hfil}}}}
\setbox\noanchorbox=\noanchor

\def\fg#1#2#3{\setbox\tempbox=\hbox{$\scriptstyle{#2}$}
\ifnum\wd\anchorbox>\wd\tempbox\dimen255=\wd\anchorbox
\else\dimen255=\wd\tempbox\fi
{#1\,\vtop{\hbox to \dimen255{\hfil\anchor\hfil}
           \vskip -6truept
           \hbox to \dimen255{\hfil$\scriptstyle{#2}$\hfil}}
           \,#3}}

\def\nfg#1#2#3{\setbox\tempbox=\hbox{$\scriptstyle{#2}$}
\ifnum\wd\noanchorbox>\wd\tempbox\dimen255=\wd\noanchorbox
\else\dimen255=\wd\tempbox\fi
{#1\,\vtop{\hbox to \dimen255{\hfil\noanchor\hfil}
           \vskip -6truept
           \hbox to \dimen255{\hfil$\scriptstyle{#2}$\hfil}}
           \,#3}}

\setbox1=\hbox{$\bot$}

\def\north#1#2{#1\,
\hbox{$\bot$\llap {\hbox to\wd1 {\hfil $/$\hfil}}}
\,#2}

\def\nao#1#2#3{#1\  \hbox{\vtop{ 
\baselineskip=4pt
\hbox{$\bot$\llap {\hbox to\wd1 {\hfil $/$\hfil}}
\hskip .05em \llap{\hbox{$^{\scriptscriptstyle{a}}$}}}\hbox{$\scriptstyle
{#2}$}}}\, #3}

\def\bp{\par{\bf Proof.}$\ \ $}

\def\includeE#1{{\lhook\kern-3.5pt\joinrel\smash{
    \mathop{\longrightarrow}\limits^{#1}}}}

\def\efor/{Example~\ref{E4}}

\def\BL/{Baldwin--Lachlan}
\def\Bu/{Buechler}
\def\Hr/{Hrushovski}
\def\lm/{locally modular}
\def\wm/{weakly minimal}
\def\nm/{non--modular}
\def\ss/{superstable}
\def\ud/{unidimensional}
\def\sm/{strongly minimal}

\def\abar{\overline{a}}

\def\bbar{\overline{b}}
\def\cbar{\overline{c}}
\def\dbar{\overline{d}}

\def\hbar{\overline{h}}

\def\mbar{\overline{m}}

\def\xbar{\overline{x}}
\def\ybar{\overline{y}}

\def\Cb{{\rm Cb}}

\def\acl{{\rm acl}}
\def\cl{{\rm cl}}

\def\tp{{\rm tp}}
\def\stp{{\rm stp}}

\def\tr/{trivial}
\def\nt/{non--trivial}
\def\st/{strong type}

\def\TV/{Tarski--Vaught}

\def\sc/{sound construction}
\def\ac/{atomic construction}

\def\fal/{functional}
\def\upl/{unique parallel lines}
\def\chp/{categorical in a higher power}

\def\abar{\bar{a}}
\def\bbar{\bar{b}}
\def\cbar{\bar{c}}
\def\dbar{\bar{d}}

\def\ybar{\bar{y}}
\def\Cbar{\overline{C}}
\def\phi{\varphi}

\def\C{{\frak  C}}

\def\FF{{\bf F}}
\def\M{{\cal M}}

\def\tp{{\rm tp}}
\def\stp{{\rm stp}}

\def\acl{{\rm acl}}
\def\dcl{{\rm dcl}}

\def\Fa0{{\FF^a_{\aleph_0}}}

\def\bp{{\bf Proof.}\quad}
\def\endproof{\hfill $\Box$ \medskip}
\def\<{\langle}
\def\>{\rangle}

\def\clp{{\rm cl}_p}
\def\Clp{\clp}
\def\Cl{{\rm cl}}
\def\clf{{\rm cl_{{\rm fork}}}}

\def\na{{\subseteq_{{\rm na}}\,}}
\newcommand\myrestriction{\mathord\restriction}
\def\mr#1{\myrestriction_{#1}}
\def\Z{{\mathbb Z}}

\newtheorem{Theorem}{Theorem}[section]
\newtheorem{Proposition}[Theorem]{Proposition}
\newtheorem{Definition}[Theorem]{Definition}

\newtheorem{Remark}[Theorem]{Remark}
\newtheorem{Example}[Theorem]{Example}
\newtheorem{Lemma}[Theorem]{Lemma}
\newtheorem{Corollary}[Theorem]{Corollary}

\newtheorem{Criterion}[Theorem]{Criterion}
\newtheorem{Fact}[Theorem]{Fact}

\begin{document}

\date{\today}

\maketitle

\abstract{ We introduce the notion of strong $p$-semiregularity and show that in a classifiable theory,  if $p$ is a regular type that is not locally modular then any $p$-semiregular type is strongly $p$-semiregular. Moreover, for any such $p$-semiregular type, ``domination implies isolation" which allows us to prove the following: Suppose that $T$ is countable, classifiable and $M$ is any model. If $p\in S(M)$ is regular but not locally modular and $b$ is any realization of $p$ then every model $N$ containing $M$ that is dominated by $b$ over $M$ is both constructible and minimal over $Mb$.}

\section{Introduction}  

In obtaining the uncountable spectrum of any classifiable theory $T$ in \cite{HHL}, localizations of $\omega$-stability near certain regular types were considered.
A regular type $p\in S(M)$ over a countable $M$ is {\em locally totally transcendental} (locally t.t.) if it is not orthogonal to a $q\in S(M)$ that is strongly regular and for which there is a constructible (and hence prime)  model over  $M$ and any realization of $q$.   There are examples of depth zero non-trivial regular types in classifiable theories which are not locally t.t (see for instance Example \ref{ex2.3}).  We intend to consider the manner in which models dominated by such types are constructed in future papers.
In this paper, we concentrate on {\em non-locally modular} regular types $p$ and prove that they are all locally t.t.\ in a very strong way.  The two main results build on
the dichotomy theorem of Hrushovski and Shelah in \cite{HrSh}.  Here, we prove that if a stationary $q\in S(A)$ is $p$-semiregular then
\begin{itemize}
\item  $q$ is strongly $p$-semiregular (see Definition~\ref{stronglysr}); and
\item  $q$ is depth-zero like and ``domination implies isolation" (DI), (see Definitions~\ref{depth0like} and \ref{DIdef}) hence if $q$ is based on a model $M$ then there is a constructible,
minimal model $N\supseteq Mb$ over $Mb$ for any realization $b$ of $q|M$,
and moreover, any $N\supseteq Mb$ that is dominated by $b$ over $M$ is constructible and minimal over $Mb$.
\end{itemize}
It is this last result which produces in particular the following easy to state theorem: 

\medskip
\noindent{\bf Theorem~\ref{atomic2}}
{\em Suppose that $M$ is any model of a classifiable theory.
If $p\in S(M)$ is regular but not locally modular and $b$ is any realization of $p$, 
then every model $N$  containing $M$ that is dominated by $b$ over $M$ is both constructible and minimal over $Mb$.}

%Although this will not be used in this paper, we note the following amusing corollary.
%
%\begin{Corollary}  \label{spectrum}  If $T$ is any countable, complete theory 
% and its spectrum of uncountable models is either A,B,C, then an infinite group is interpretable in models of $T$.
%In particular, $I(T,\aleph_{\alpha})=\min\{2^{\aleph_\alpha},\beth_2\}$ (which e.g., is the spectrum of $Th(\Z,+)$), then $T$ must interpret an infinite group.
%\end{Corollary}
%
%\bp  BY XYZ, if $T$ does not interpret an infinite group, then every non-trivial regular type cannot be locally modular.   Thus, by Theorem~\ref{atomic2}, 
%any non-trivial regular type $p\in S(M)$ is locally t.t.\ in $M$.  So, in the notation of \cite{1}, for any $n\ge 1$, 
%`non-trivial failure' $NTF(n)$ is impossible, hence $TT(n)$ holds if and only if $TF(n)$ fails.
%The result follows by considering how the spectra for A,B,C are attained in Section~XXXXXXXX

\medskip

We assume that the reader is familiar with superstability and adopt the usual convention when working with stable theories that we are working in a large, saturated model $\C$ of the theory and all models mentioned are small elementary submodels of $\C$, sets are small subsets of $\C$ and tuples are from $\C$. We will also assume that $T$ eliminates imaginaries i.e. $T = T^{eq}$.
This paper has an extensive Appendix that records a number of  definitions and facts from basic geometric stability theory.  
Section \ref{Prelim} contains some background information and gives an example indicating the necessity of non-local modularity in Theorem~\ref{atomic2}.  Section 3  introduces
the notion of strongly $p$-semiregular strong types and proves that when $p$ is any non-locally modular regular type, every $p$-semiregular  strong type is strongly $p$-semiregular.
We  provide some applications of this in Section 4 before proving our main theorem in Section 5.   The paper proper concludes with an examination of other circumstances under which  domination implies isolation in Section 6.

We are grateful to the anonymous referee, as his/her thorough reading of this paper has led to many expositional improvements.

\medskip
\centerline{{\bf For the whole of this article we assume $T$ is (at least) stable.}}

\medskip

\section{A brief historical background}\label{Prelim}

One of the major accomplishments of stability theory was Shelah's proof of the Main Gap in \cite{Shc} where the notion of a classifiable theory was introduced.  

\begin{Definition} {\em  A complete theory $T$ in a countable language is {\em classifiable} if $T$ is superstable, has prime models over pairs (PMOP) and does not have the dimension order property (NDOP) 
}
\end{Definition}

The definitions of these terms, along with equivalents we will use are given in Definition~\ref{basicdef} and Facts \ref{NDOP} and \ref{PMOPequiv}.
The significance of classifiability is given by two results.  On one hand, in \cite{ShelahBuechler}, building extensively on \cite{Shc}, Shelah and Buechler prove
that if $T$ is classifiable, then every model of $T$
 is constructible\footnote{See Section~\ref{appendix-domisol} of the Appendix.}
and minimal over
an  independent tree of countable elementary substructures. 
By contrast,  
if a countable theory is not classifiable, then $I(T,\kappa)=2^\kappa$ for every uncountable cardinal $\kappa$.  The reader is cautioned that these two properties are not exclusive
--
a classifiable, deep theory has both  the structure theorem and has
 maximally many non-isomorphic models in every uncountable cardinality.

%The two properties, NDOP and PMOP, refer to independent triples $\M=(M_0,M_1,M_2)$ of models of $T$,
%where $M_0\subseteq M_1$, $M_0\subseteq M_2$, and $\fg {M_1} {M_0} {M_2}$.  A superstable theory {\em does not have the dimension order property} (NDOP) if, for every triple $\M$ of a-models,
%the a-prime model $M^*$ over $M_1M_2$ (which exists in any superstable theory) is minimal among all a-models containing $M_1M_2$.
%Shelah proves that NDOP is equivalent to the statement that any regular type $q$ non-orthogonal to $M^*$ is either non-orthogonal
%to $M_1$ or $M_2$.  
%
%A superstable theory $T$ has {\em prime models over pairs } (PMOP)  if, for any independent triple $(M_0,M_1,M_2)$ of models, there is a 
%constructible\footnote{Recall that a model $N$ is constructible over a set $B$ if its universe can be enumerated as $\{c_i:i<\alpha\}$ with $\tp(c_i/B\cup\{c_j:j<i\})$ isolated for each $i$.
%Any two constructible models over $B$ are isomorphic.}
% model
%over $M_1M_2$.  

After Shelah defined PMOP, Harrington gave an alternate treatment of PMOP that  was developed in \cite{Hart}.  Following the notation preceding Fact~\ref{PMOPequiv},
call a type $p\in S(A_1A_2)$ over an independent triple
$(A_0,A_1,A_2)$ of sets {\em $V$-dominated} if, for every extension $(B_0,B_1,B_2)$ of $(A_0,A_1,A_2)$, and for any realization $c$ of $p$ with $\fg c {A_1A_2} {B_0}$,
we have $\fg c {A_1A_2} {B_1B_2}$.   Then, as in Fact~\ref{PMOPequiv}(4),  Harrington notes that the statement ``Every $V$-dominated type is isolated" is equivalent to PMOP for countable, superstable theories with NDOP.  
The form of this statement is indicative of what we do here.  In Section 6 we investigate conditions under which `domination implies isolation' in various settings.

%Given two independent triples $\M=(M_0,M_1,M_2)$, $\N=(N_0,N_1,N_2)$
%of models, say that $\N$ {\em extends $\M$} if $M_0\subseteq N_0$, $\fg {N_0} {M_0} {M_1M_2}$, and $N_i$ is dominated by $M_i$ over $N_0$ for $i=1,2$.
%Call a  type $\tp(b/M_1M_2)$ {\em $V$-isolated} if $\fg b {M_1M_2} {N_1N_2}$ for every extension $\N$ of $\M$.
%Among countable superstable theories with NDOP, PMOP is equivalent to ``every $V$-isolated strong type is isolated."
%
%This  notion of $V$-isolation  will not be used directly in the paper, but it is  an example of how, in classifiable theories, certain forms of domination or ``weak'' isolation,  imply actual  isolation. This is a theme that occurs throughout  the paper and is highlighted in Section 6.

%In \cite{ShelahBuechler} but building extensively on \cite{Shc},
%Shelah and Buechler prove that among complete countable theories $T$, $T$ is classifiable if and only if every model is prime and minimal over an
%independent tree of countable, elementary substructures.  We shall call such a tree a classifying tree.

Since the proof of the structure theorem, there has been a considerable amount of work analyzing the `fine structure' of classifiable theories.  The fine structure to a large extent revolves around understanding the leaves of classifying trees.  The leaves are controlled by depth zero types and so we remind the reader of the definition (for more definitions and classical results consult the Appendix). 
\begin{Definition} \label{depth0def}
{\em A regular type $p$ in a superstable theory is said to have {\em depth zero} if for any a-model $M$ on which $p$ is based and any realization $b$ of $p|M$, any non-algebraic type $q$ over $M[b]$ (the a-prime model over $Mb$), is non-orthogonal to $M$ i.e. $p$ does not support a regular type. }
\end{Definition}
A {\it leaf} is a triple $(M,b,N)$ where $M \preceq N$, $b \in N$, $tp(b/M)$ is regular and depth zero, and $N$ is dominated by $b$ over $M$. The computation of the uncountable spectrum of a countable theory depends on, among other things, understanding the isomorphism types of leaves that appear in the classifying trees of models.  Through the coarseness of cardinal arithmetic, in \cite{HHL}, it was not necessary to determine all possible isomorphism types of leaves in order to determine the uncountable spectra of countable theories.  Still, the techniques available from classification theory are suitable to do this and demonstrate a deeper understanding of the structure of models of classifiable theories.  In this paper and subsequent papers, we intend to explore the isomorphism types of leaves $(M,b,N)$ at least when $tp(b/M)$ is non-trivial; in a classifiable theory, non-trivial types necessarily have depth zero (Fact~\ref{NDOP}).   Our main result here is that if $(M,b,N)$ is a leaf with $\tp(b/M)$ regular and non-locally modular, then
there is only one choice for the isomorphism type of $N$ over $Mb$.    This leads to an amusing corollary to Theorem~\ref{atomic2}.

%Although this will not be used in this paper, we note the following amusing corollary.

\begin{Corollary}  \label{spectrum}  If $T$ is any countable, complete theory 
 and its spectrum of uncountable models is any of 4,5,7, or 10 (in the notation of Theorem~6.1 of \cite{HHL}) then an infinite group is interpretable in models of $T$.
In particular, if $I(T,\aleph_{\alpha})=\min\{2^{\aleph_\alpha},\beth_2\}$ (i.e., the spectrum of $Th(\Z,+)$), then $T$ must interpret an infinite group.
\end{Corollary}

\bp   Hrushovski proves that  if $T$ does not interpret an infinite group, then every non-trivial regular type is non-locally modular, see e.g., 7.3.3 of \cite{Pillay}.
 Thus, by 
Theorem~1b of \cite{HrSh} and Theorem~\ref{atomic2} here, 
any non-trivial regular type $p\in S(M)$ is locally t.t.\ in $M$.  So, in the notation of \cite{HHL}, for any $n\ge 1$, 
`non-trivial failure' $NTF(n)$ is impossible, hence $TT(n)$ holds if and only if $TF(n)$ fails.
The result follows by considering the case analysis for each spectrum given in Section~5 of \cite{HHL}.   
\endproof

%We recall two separate facts about leaves that were known before this investigation.
%
%First of all, if $(M,b,N)$ is a leaf in a classifiable theory, $tp(b/M)$ is non-trivial and $M \subseteq_a \C$ (all strong types of finite tuples from $\C$ over finite tuples in $M$ are realized in $M$) then the isomorphism type of $N$ is determined up to isomorphism over $M$; see \cite{Shc}.  Of course such an $M$ would typically have size at least the continuum but this still shows that the geometry of $tp(b/M)$ plays a role in the isomorphism types of leaves.
%
%Secondly, recall that if $T$ is countable and stable, then for any set $A$, there is an $\ell$-constructible model $N$ over $A$ (see Definition in Section \ref{defisolation}).  In particular, if $M \subset N$ and $N$ is $\ell$-constructible over
%$Mb$, then $N$ is dominated by $b$ over $M$.    In fact, if $T$ is countable, superstable, and NDOP, and $(M,b,N)$ is a leaf then $N$ is $\ell$-constructible over $Mb$. This doesn't say that the isomorphism type of $N$ is determined by $Mb$ but it does put constraints on how such $N$ can be built.

%Our goal is to measure the extent to which these results can be extended to constructible models assuming classifiability.
The results of this paper are concerned with non-locally modular regular types over models.  It is natural to ask for other circumstances under which leaves are constructible. This is always possible if $T$ is $\omega$-stable, as such a theory has constructible models over any set.
However, we close this section with an example of a classifiable theory that has leaves $(M,b,N)$ for which $N$ is not constructible over $Mb$.
%Within the context of $\omega$-stable theories, this is easy.  
% As $\omega$-stable theories have constructible models over any set,
%an $\omega$-stable theory is classifiable if and only if it has NDOP.  For such theories, if $(M,b,N)$ is a leaf then $N$ is constructible over $Mb$.  Unfortunately, for an arbitrary classifiable theory, there are leaves $(M,b,N)$ for which $N$ is not constructible over $Mb$, as the following example shows.

\begin{Example} \label{ex2.3} {\em The language will consist of countably many sorts $\{U_n:n\ge 1\}$,   a collection of unary relations $\{R^n_\eta:\eta \in 2^{<\omega},n\ge 1\}$,  and on each sort $U_n$, there is a binary operation $+_n:U_n^2\rightarrow U_n$.  The language also includes projection functions $f_n:U_{n+1}\rightarrow U_{n}$
for each $n\ge 1$. The canonical model of our theory in this language is as follows:
\begin{enumerate}
\item $U_n$ is interpreted as the product of $n$ copies of $2^\omega$;
\item $+_n$ is interpreted as coordinate-wise addition modulo 2;
\item If $b\in U_n$ codes the  $n$-tuple $\nu_1\ldots,\nu_{n}$, then $R^n_\eta(b)$ holds  iff $\eta$ is an initial segment of $\nu_n$, and
\item $f_n$ is the projection onto the first $n$ coordinates.
\end{enumerate}
The theory $T$ of this structure is classifiable; in fact, it is superstable and unidimensional but not $\omega$-stable.  Now suppose that $M$ is a model of $T$ and $N$ is an elementary weight one extension of $M$ - the type of thing that would happen with leaves on a classifying tree.  The claim is that if $b \in N \setminus M$ is any finite tuple then $N$ is not constructible over $Mb$.  To see this, suppose we have such a $b$.  By the presence of addition in all the sorts and the model $M$, we can assume that $b$ is a singleton in some sort $U_n$.  But then, if one considers the preimage of $b$ under $f_{n}$, one sees that this formula in the sort $U_{n+1}$ does not contain an isolated type - the predicates $R^{n+1}_\eta$ preclude this.}
\end{Example}

This example is suggestive of the result we will prove in subsequent papers: if we don't restrict ourselves to finite tuples then of course $N$ is determined by making a coordinated choice of elements from each sort.  
In the example, all the regular types are locally modular (non-trivial).   It is not clear in advance that this is important but in this paper we will show that if $(M,b,N)$ is a leaf in a countable, classifiable theory and $tp(b/M)$ is {\em not}  locally modular, then $N$ is constructible over $Mb$.

%
%\medskip
%\centerline{{\bf For the whole of this paper  we assume $T$ is (at least) stable.}}
%\medskip
%

\section{Strongly $p$-semiregular types}

Definitions of  $p$-semiregular types,  $p$-simplicity and other related notions can be found in the Appendix.

\begin{Definition} \label{stronglysr}  {\em  A stationary type $q\in S(B)$ is {\em strongly $p$-semiregular} of weight $k>0$  if
$q$ is $p$-semiregular of $p$-weight $k$  and there is a $p$-simple formula $\theta(x)\in q$  of weight $k$ such that
if $d$ realizes $\theta$ and $C\supseteq B$ with $\fg d B C$ and $w_p(d/C)=w_p(q)$, then $\tp(d/C)=q|C$, the non-forking 
extension of $q$ to $S(C)$.
}
\end{Definition}

The goal of this entire section is to prove the following Theorem.  Its proof is patterned after the argument in \cite{HrSh}, where Hrushovski and
Shelah prove that in a classifiable theory, every non-locally modular stationary, regular type is strongly regular (i.e., strongly $p$-semiregular of
$p$-weight one). 
% We will refer to precise Lemmas in their paper quite a few times in this section. When there is a risk of ambiguity,  these Lemmas will be denoted as ``Lemma~x.y\cite{HrSh}''.

\begin{Theorem}  \label{strong}  If $T$ is classifiable and  $p$ is a non-locally modular regular type, then every  stationary $p$-semiregular type is
strongly $p$-semiregular.
\end{Theorem}

%\begin{Remark}  {\em  It appears that all is needed is $T$ superstable, PMOP, and $p$ is non-locally modular of depth zero.
%}
%\end{Remark}

\begin{Remark} \label{add}  {\em  By the open mapping theorem, this notion is parallelism invariant.  In particular, if $d \anchor E$
with $\tp(d/\emptyset)$ stationary, then if $\tp(d/E)$ is strongly $p$-semiregular via $\theta(x,e)$, then $\tp(d/\emptyset)$
will be strongly $p$-semiregular via $\bigvee d_r y\theta(x,y)$, the (finite) disjunction of $\theta$-definitions of each $r\in S_\theta(\C)$ consistent with $\tp(e/\emptyset)$.
}
\end{Remark}

We will use this Remark as justification for freely adding independent parameters in many places.

\subsection{$p$-triples}  In this subsection, $T$ is superstable and $p$ is a stationary regular type over $\emptyset$.  The reader is encouraged to review the definitions and basic facts in Sections~\ref{appendix-locallymodular} and \ref{appendix-semiregular} of the Appendix.
We adopt the data structure of a {\em $p$-triple} and then show that a given $p$-triple can be massaged to get matching  $p$-triples with more and more
desirable properties.  The key will be to obtain a {\em minimal $p$-triple} as a normal cover of a given one.

\begin{Definition}
{\em  
\begin{itemize}
\item  A {\em $p$-triple} is a sequence $(a,b,C)$ such that  $a,b\in D(p,C)$, i.e., $\stp(ab/C)$ is $p$-simple, with both $w_p(a/Cb)>0$ and $w_p(b/Ca)>0$.
\item  A $p$-triple $(a,b,C)$ is {\em normal} if the three strong types $\stp(ab/C)$, $\stp(a/Cb)$, and $\stp(b/Ca)$ are all $p$-semiregular.
\item  A $p$-triple $(a,b,C)$ is {\em $p$-disjoint} if $\Clp(Ca)\cap \Clp(Cb)=\Clp(C)$.
\item  Two $p$-triples $(a,b,C)$, $(a',b',C')$ are {\em matching} if $w_p(ab/C)=w_p(a'b'/C')$, $w_p(a/Cb)=w_p(a'/C'b')$, and $w_p(b/Ca)=w_p(b'/C'a')$.
\end{itemize}
}
\end{Definition}

We seek $p$-triples $(a,b,C)$ that embed a realization of a given $p$-semiregular strong type in its first component.

\begin{Definition}  {\em  Suppose $\stp(e/\emptyset)$ is $p$-semiregular.  
A $p$-triple $(a,b,C)$  {\em envelops $e$} if
$e\in\dcl(a)$ and $e \anchor {Cb}$.
}
\end{Definition}

There is a canonical way of extending a $p$-disjoint $p$-triple $(a,b,C)$ to a matching, normal $p$-disjoint $p$-triple $(a',b',C')$. 

\begin{Definition} {\em  A  {\em soft extension} of a $p$-triple $(a,b,C)$ is a $p$-triple $(a',b',C')$ such that 
\begin{enumerate}
\item  $a'b'C'\subseteq\dcl(abC)$;
\item  $C\subseteq C'\subseteq\clp(C)$  with $C'\setminus C$ finite; 
\item  $a\subseteq a'\subseteq \Clp(C'a)$ and
 $b\subseteq b'\subseteq \Cl_p(C'b)$.
\end{enumerate}
}
\end{Definition}

Note that $\clp(C')=\clp(C)$, $\clp(C'a')=\clp(Ca)$, and $\clp(C'b')=\clp(Cb)$ in any soft extension.

\begin{Lemma}[Normalization]\label{normalization}
Given any $p$- triple $(a,b,C)$, there is a soft extension to a  matching normal $p$-triple $(a',b',C')$, called a normalization of $(a,b,C)$.
If $(a,b,C)$ is $p$-disjoint, then $(a',b',C')$ will be $p$-disjoint as well.   Moreover, if $\stp(e/\emptyset)$ is $p$-semiregular and $(a,b,C)$ envelops $e$,
then $(a',b',C')$ also envelops $e$.

\end{Lemma}

\bp  First, apply Lemma~\ref{sr-existence} to $\stp(ab/C)$ to get $C'\subseteq\dcl(Cab)\cap\clp(C)$
such that $\stp(ab/C')$ is $p$-semiregular.  Note that $\Clp(C')=\Clp(C)$,
$\Clp(C'a)=\Clp(Ca)$, and $\Clp(C'b)=\Clp(Cb)$. Next, apply the Lemma to
$\stp(b/aC')$ to get $a'$ and $\stp(b/C'a')$ $p$-semiregular, and finally 
apply the Lemma to $\stp(a'/C'b)$ to get $b'$.  It follows from Lemma~\ref{semiregpersists}  that $(a',b',C')$ is normal and matches $(a,b,C)$.
The preservation of $p$-disjointness is clear because of the equality of the $p$-closed sets mentioned above.
Finally, suppose $(a,b,C)$ envelops $e$.  Since $a\subseteq a'$, $e\in\dcl(a')$.  Also, since $C'b'\subseteq\Clp(Cb)$, $\stp(e/\emptyset)$ $p$-semiregular and
$e \anchor {Cb}$ implies $ e \anchor {C'b'}$, so $(a',b',C')$ envelops $e$ as well.
\endproof

Next we describe three ways of extending a given $p$-triple $(a,b,C)$ to a larger, matching $(a',b',C')$ that preserves $p$-disjointness and enveloping.

\begin{Definition}  {\em  A {\em simple extension} of a $p$-triple $(a,b,C)$ is any of:
\begin{enumerate}
\item $(a^*,b,C)$, where $a\subseteq a^*\subseteq \Clp(Ca)$;
\item $(a,b^*,C)$, where $b\subseteq b^*\subseteq \Clp(Cb)$;
\item  $(a,b,C^*)$, where $C^*\supseteq C$ and $\fg {C^*} C {ab}$.

\end{enumerate}

}
\end{Definition}

\begin{Lemma}  \label{simple}  Given any $p$-disjoint, normal  $p$-triple $(a,b,C)$, the normalization 
$(a',b',C')$ of any simple extension is a matching $p$-disjoint extension of $(a,b,C)$.
Moreover, if $\stp(e/\emptyset)$ is $p$-semiregular and $(a,b,C)$ envelops $e$,
then $(a',b',C')$ envelops $e$ as well.

\end{Lemma}

\bp  That all three species of simple extensions are matching is clear.  Next, we show that each of the simple extensions preserves $p$-disjointness.
This is clear for the first two, as $\clp(Ca^*)=\clp(Ca)$ in the first case and $\clp(Cb^*)=\clp(Cb)$ in the second.  As $(a,b,C)$ normal implies $\stp(ab/C)$
is $p$-semiregular, $p$-disjointness of the third species is preserved by Lemma~\ref{nonforking}.  If $(a,b,C)$ envelops $e$, then obviously $e\in\dcl(a^*)$ for any of the three species of simple extensions.  For the first species, $ e \anchor {Cb}$ holds vacuously.  For the second we have $Cb^*\subseteq\Clp(Cb)$, so $e \anchor {Cb^*}$
follows from $\stp(e/\emptyset)$ $p$-semiregular, and for the third, since $\fg {eb} {C} {C^*}$, we obtain $ e \anchor  {C^*b}$ by the transitivity of non-forking.

The Lemma now follows from Lemma~\ref{normalization}.
\endproof

\begin{Definition} 
{\em  Suppose that $(a,b,C)$ is a $p$-disjoint normal $p$-triple.    
\begin{itemize}
\item  A {\em normal cover} is any  $p$-disjoint normal $(a',b',C')$ obtained as a sequence of extensions as in Lemma~\ref{simple}.
\item  The {\em strength}  of $(a,b,C)$ which we denote by $\alpha(a,b,C)$, is equal to $R^\infty(a/bb'C)$,
where $b'$ is (any) element satisfying $\stp(b'/Ca)=\stp(b/Ca)$ and
$\fg {b'} {Ca} b$.  (This is well-defined, as $\stp(abb'/C)$ is independent of
our choice of $b'$.)

\item   A normal $p$-triple $(a,b,C)$ is {\em minimal\/}
if $\alpha(a,b,C)\le\alpha(a',b',C')$ for all  of its normal covers $(a',b',C')$.
\end{itemize}
}
\end{Definition}

Clearly, by superstability and the transitivity of being a normal cover, every normal $p$-triple $(a,b,C)$ has a minimal, normal cover 
$(a',b',C')$.  
In fact, one can find one with the additional property
 that $a'\in\dcl(C'a)$.  At present, this improvement does not seem to be necessary.

\begin{Lemma}  \label{minimalpass}  Given any  $p$-disjoint, normal
$p$-triple $(a,b,C)$, there is a  matching minimal, normal, $p$-disjoint $p$-triple $(a',b',C')$ that is a normal cover of $(a,b,C)$.
If $\stp(e/\emptyset)$ is $p$-semiregular and $(a,b,C)$ envelops $e$, then $(a',b',C')$ envelops $e$ as well.
\end{Lemma}

\bp  Among all normal covers $(a',b',C')$ of $(a,b,C)$, choose the one of smallest strength.  \endproof

We record how the minimality assumption will be used in the proof of Theorem~\ref{strong}.  

\begin{Lemma}  \label{howused}  Suppose $(a,b,C)$ is a minimal normal $p$-triple with strength $\alpha$.
For every $\Cbar\supseteq C$ with $\fg {\Cbar} C {ab}$ and every $B$ with $b\subseteq B\subseteq \clp(\Cbar b)$,
$R^\infty(a/\Cbar B B')=\alpha$ for some/every $B'$ with $\stp(B'/\Cbar a)=\stp(B/\Cbar a)$ and $\fg B {\Cbar a} {B'}$.
\end{Lemma}

\bp  This is immediate since $(a,b,\Cbar)$ is a simple extension of  $(a,b,C)$, and $(a,B,\Cbar)$ is a simple extension of $(a,b,\Cbar)$, hence 
$(a,B,\Cbar)$ is a normal cover of $(a,b,C)$.
\endproof

We close with two lemmas concerning $p$-disjointness.  

\begin{Lemma}  \label{pindep}
Suppose $\stp(a/C)$ and $\stp(b/C)$ are both $p$-semiregular and $\fg a C b$.
Then the $p$-triple $(a,b,C)$ is  $p$-disjoint.
\end{Lemma}

\bp    Choose any $e\in\Clp(Ca)\cap\Clp(Cb)$.
Then $w_p(e/Ca)=w_p(e/Cb)=0$.  By Lemma~\ref{semiregpersists}(1), $\stp(b/Ca)$ is $p$-semiregular,
so $\fg b C e$.  But this, coupled with 
$w_p(e/Cb)=0$ implies $e\in\Clp(C)$.
\endproof

\begin{Lemma}  \label{pair}  Suppose $(a_1,b_1,C)$ and $(a_2,b_2,C)$
are both  normal and $p$-disjoint.
If, moreover, $\fg {a_1b_1} C {a_2b_2}$, then the $p$-triple
$(a_1a_2,b_1b_2,C)$ is also  $p$-disjoint.
\end{Lemma}

\bp  In light of Lemma~\ref{nonforking},
then $p$-disjointness of $(a_1,b_1,C)$
implies that $\Clp(a_1a_2b_2C)\cap\Clp(b_1a_2b_2C)\subseteq\Clp(a_2b_2C)$,
so $$\Clp(a_1a_2C)\cap\Clp(b_1b_2C)\subseteq\Clp(a_2b_2C)$$
Arguing in reverse, the  $p$-disjointness of $(a_2,b_2,C)$
yields
 $$\Clp(a_1a_2C)\cap\Clp(b_1b_2C)\subseteq\Clp(a_1b_1C)$$
Furthermore, it follows immediately from Lemma~\ref{pindep} that
$$\Clp(a_1b_1C)\cap\Clp(a_2b_2C)\subseteq\Clp(C)$$
and the result follows.  
\endproof

\subsection{Proof of Theorem~\ref{strong}}

We proceed to prove Theorem~\ref{strong}, following the proof of Proposition~3.1 in \cite{HrSh}.  
For the whole of this subsection, assume $T$ is classifiable and $p$ is a (stationary) regular, non-locally modular type $p$ over $\emptyset$.  By Fact~\ref{NDOP}, $p$ has depth zero.  We will often use Fact~\ref{pweightprov}, which says that every $p$-semiregular type contains a  $p$-simple formula $D$ inside which $p$-weight is continuous and definable.  
Fix some $e$ such that $\tp(e/\emptyset)$ is stationary and $p$-semiregular with $w_p(e)=k\ge 1$.  
Our goal is to show that $\tp(e/\emptyset)$ is strongly $p$-semiregular.    For the whole of this subsection, we work inside $D^{eq}$, which suffices.  Thus, all types under consideration will be $p$-simple.

%We recall in the appendix basic definitions and facts about regular locally modular types (Section \ref{appendix-locallymodular})

The following Proposition is the content of Lemma~3.2 of \cite{HrSh}.

\begin{Proposition}  \label{basic}
For every realization $d$ of $p$, there
is a minimal, normal,  $p$-disjoint $p$-triple $(a,b,C)$
enveloping $d$.   Moreover, $w_p(a/C)=w_p(b/C)=2$, $w_p(ab/C)=3$, and
$w_p(a/bC)=w_p(b/aC)=1$.
\end{Proposition}

\bp  By the first paragraph of the proof of Lemma~3.2 in \cite{HrSh},
there is a $p$-triple $(a,b,M)$, where $M$ is an a-model, $a$ and $b$ each consist
of two $M$-independent realizations of $p|M$, with $w_p(a/Mb)=w_p(b/Ma)=1$
and the$p$- triple $(a,b,M)$ being $p$-disjoint.  

By  $p$-weight computations,
at least one of $\{a_1,a_2\}$ must be independent from $b$ over $M$,
so by passing to an automorphism of $\C$, we may assume that $d\subseteq a$
and $ d \anchor {Mb}$.
Thus, $(a,b,M)$  envelops $d$.  Now, apply the Normalization Lemma~\ref{normalization}, and then choose a minimal normal cover
via Lemma~\ref{minimalpass}.
\endproof

The following generalization is the point of all these definitions.

\begin{Proposition}  \label{k}
For $e$ as above,
there is a minimal, normal, $p$-disjoint $(a,b,C)$ enveloping $e$
such that $w_p(ab/C)=3k$, $w_p(a/C)=w_p(b/C)=2k$, and $w_p(a/bC)=w_p(b/aC)=k$.
\end{Proposition}

\bp  First, by $p$-semiregularity,  choose $E$ and an $E$-independent
sequence $\dbar=(d_i:i<k)$ of realizations of $p|E$ such that $e \anchor E$
and $e$ and $\dbar$ are domination equivalent over $E$.  
%Without loss, we may assume $E=\emptyset$.  

By Proposition~\ref{basic}, choose a triple
$(a_0,b_0,C_0)$ enveloping $d_0$ as there.  Since $d_0$ realizes $p|C_0$ and $p|E$, we may assume $E\subseteq C_0$ and $\dbar$ is independent over $C_0$.
% Since $\fg d {\phantom{e}} {C_0}$, we may assume without loss of generality that
%$\{d_i:i<k\}$ are an independent set of realizations of $\stp(d_0/C_0)$.  
%Without loss, we may
%assume that it envelops $d_0$.  
Next, choose $(a_ib_i:i<k)$ to be 
$C_0$-independent with $\stp(a_ib_id_i/C_0)=\stp(a_0b_0d_0/C_0)$ for each $i$.
In particular,  $(a_i,b_i,C_0)$ envelops $d_i$ for each $i<k$.
As notation, let $\abar$ denote $(a_i:i<k)$ and $\bbar$ denote
$(b_i:i<k)$.  Note that by the independence, $\fg {b_i} {C_0a_i}  {\abar}$ for each $i$.

By iterating Lemma~\ref{pair}, we have that the $p$-triple $(\abar,\bbar,C_0)$ is
 $p$-disjoint.  As well, it follows from the independence
that $w_p(\abar/C_0)=w_p(\bbar/C_0)=2k$, $w_p(\abar\bbar/C_0)=3k$,
and $w_p(\abar/C_0\bbar)=w_p(\bbar/C_0\abar)=k$.
It also follows that $(\abar,\bbar,C_0)$ is normal.

As $w_p(e/d_0,\dots,d_{k-1})=0$, we have
$w_p(e/ \abar C_0)=0$.  Thus, the $p$-triple $(e\abar,\bbar,C_0)$ is a simple
extension of $(\abar,\bbar,C_0)$ enveloping $e$.   By  Lemmas~\ref{normalization} and \ref{minimalpass} there is a matching, minimal, normal cover
$(a,b,C)$ of $(e\abar,\bbar,C_0)$ that is  $p$-disjoint and envelops $e$.
\endproof

%We continue  along the lines of Section 3 of \cite{HrSh}, with our
%Proposition~\ref{k} taking the place of their Lemma~3.2.  A key point is
%that we have the following Lemma, which takes the place of their Lemma~3.3.

\medskip

%\begin{quotation}
{\bf From now on,  fix a $p$-triple $(a,b,C)$ as in Proposition~\ref{k}.}
%\end{quotation}

\medskip

\noindent We continue  along the lines of Section 3 of \cite{HrSh}, with our
Proposition~\ref{k} taking the place of their Lemma~3.2.  
The following Lemma  takes the place of their Lemma~3.3.

\begin{Lemma}  \label{getisol} Choose any $b'$ realizing $\stp(b/Ca)$ with $\fg {b'} {Ca} b$.
Then:
\begin{enumerate}
\item  $w_p(a/bb'C)=0$;
\item  $\fg {b'} C b$;
\item  There is a formula $\rho(x,b,b')$ (possibly with hidden parameters from $C$) isolating $\tp(a/Cbb')$.
\item  For any $a'$ and any $b''$ satisfying $\stp(b''/C)=\stp(b/C)$ and $\fg {b''} C b$, if $\rho(a',b,b'')$ holds, then $\tp(a'/Cb)=\tp(a/Cb)$.
\end{enumerate}
\end{Lemma}

\bp  (1)  Here is where $p$-disjointness plays a leading role.
Recalling the notation from the proof of Proposition~\ref{k}, 
$(a,b,C)$ is a normal cover of $(e\abar,\bbar,C_0)$,  and for each
$i<k$ we have $w_p(a_i/b_iC_0)=1$ and $\fg {b_i} {C_0a_i}  {\abar}$.
We claim that $\nfg {a_i} {C_0b_i} {b'}$ for each $i<k$.  Since $\stp(a_i/b_iC_0)$ is $p$-semiregular, the claim implies that $w_p(a_i/b_ib'C)=0$ for each $i$, hence
$w_p(a/bb'C)=0$, completing the verification of (1).

To establish the claim,  fix $i<k$.  The forking is attained in three steps.
%
%
%As $\stp(a_i/b_iC_0)$ is $p$-semiregular,
% we additionally have $w_p(a_i/bC)=1$ as well.
%So, in order to establish (1), it suffices to prove that $w_p(a_i/bb'C)=0$
%for each $i<k$.
%
%Fix an $i<k$.  We will actually prove that $\nfg {a_i} {C_0b_i} {b'}$, which suffices.  To see this, 
First, by $p$-disjointness we have
$$\Clp(Ca)\cap\Clp(Cb)\subseteq\Clp(C)$$
Thus, $\acl(C_0a_i)\cap\acl(C_0b_i)\subseteq\Clp(C)$.  Coupling this with the fact
that $\stp(a_ib_i/C_0)$ is $p$-semiregular implies $\fg {a_ib_i} {C_0} {\Clp(C)}$, hence
$$\acl(a_iC_0)\cap \acl(b_iC_0)\subseteq\acl(C_0)$$
Second, from the proof of Proposition~\ref{k}, $\fg {b_i} {C_0a_i} {\abar}$.  By the normality of $(a_i,b_i,C_0)$, $\stp(b_i/C_0a_i)$ is $p$-semiregular,
so since $Ca\subseteq\Clp(C_0\abar)$ we have $\fg  {b_i} {Ca_i} a$.  Our assumption was that $\fg b {Ca} {b'}$, so in particular, $\fg {b_i} {Ca} {b'}$.
Thus, by the transitivity of non-forking
$$\fg {b'} {Ca_i} {b_i}$$
Finally, since $\stp(b'/Ca)=\stp(b/Ca)$, there is $b_i'\in b'$
(corresponding to $b_i$) such that $\nfg {b_i'} {C_0} {a_i}$,
hence
$$\nfg {b'} {C_0} {a_ib_i}$$
Combining the last three displayed expressions with Lemma~\ref{fork}
(where $b'$ takes the role of $X$) we obtain $\nfg {a_i} {C_0b_i} {b'}$.

(2)  Given (1), this  is analogous to the $p$-weight computation given in Lemma~3.3(b) of \cite{HrSh}, just multiplied by $k$.  First, 
 $w_p(abb'/C)=w_p(a/C)+2w_p(b/Ca)=4k$.  Using (1),  $4k=w_p(abb'/C)=w_p(bb'/C)$. Since $w_p(b/C)=w_p(b'/C)=2k$, we conclude
 that $w_p(bb'/C)=w_p(b/C)+w_p(b'/C)$.    As $\stp(b/C)$ is $p$-semiregular, this implies
 $\fg b C {b'}$.

(3)  This is analogous to 3.3(c) of \cite{HrSh}. We employ Fact~\ref{PMOPequiv}(4) and the notation immediately preceding it.
By (2), $(C,Cb,Cb')$ is an independent triple, so  it suffices to show that $\tp(a/Cbb')$ is $V$-dominated.
Choose any independent triple $(\Cbar,B,B')$ extending $(C,b,b')$, i.e., 
 $\Cbar\supseteq C$ with $\fg {\Cbar} C {abb'}$ and $B\supseteq \Cbar b$, $B'\supseteq \Cbar b'$ with $\fg B  {\Cbar} {B'}$.  
We show that $\fg a {\Cbar b b'} {BB'}$ by partitioning $BB'$.
into its $p$-weight zero part and its $p$-semiregular part.  
More precisely, let  $B_0:=\dcl(B)\cap\Clp(\Cbar b)$ and $B_0':=\dcl(B')\cap\Clp(\Cbar b')$.  
Since $\fg b {\Cbar a} {b'}$, there is an automorphism $\sigma\in Aut(\C)$ fixing $\acl(\Cbar a)$ pointwise and $\sigma(b)=b'$.  Note that $\sigma$ maps $\Clp(\Cbar b)$ onto $\Clp(\Cbar b')$, so we can find $B_1$ with $B_0\subseteq B_1\subseteq \Clp(\Cbar b)$ and $B_0'\subseteq \sigma(B_1)$.
As our original normal triple $(a,b,C)$ is minimal, it follows from Lemma~\ref{howused} that $R^\infty(a/\Cbar B_1\sigma(B_1))=R^\infty(a/C bb')$.  In particular,
$$\fg a {\Cbar b b'} {B_0B_0'}$$

Continuing, recall that we are working entirely within a $p$-simple formula $D$.  Thus, by Lemma~\ref{sr-existence}, both $\stp(B/B_0)$ and $\stp(B'/B_0')$ are $p$-semiregular,
with at least one of positive $p$-weight.  
As well,  $\Cbar\subseteq B_0\cap B_0'$ and $\fg {B} {\Cbar} {B'}$, hence $\stp(BB'/B_0B_0')$ is $p$-semiregular of positive $p$-weight as well.  
By (1)  $w_p(a/B_0B_0')=0$, so $\fg a {B_0B_0'} {B_1B_1'}$.  Thus,  $\stp(a/Cbb')$ is $V$-dominated, and hence by Fact~\ref{PMOPequiv}(4)
there is $\rho(x,b,b')$ (possibly with hidden parameters from $C$) isolating $\tp(a/Cbb')$.

(4)  By (2), any such $b''$ satisfies $\tp(bb''/C)=\tp(bb'/C)$, so $\rho(a',b,b'')$ implies $\tp(a'/Cb)=\tp(a/Cb)$.
\endproof

%We record the following consequence of Lemma~\ref{getisol}.
%
%\begin{Conclusion}  \label{whatweuse}  For any $a'$ and  any $b''$ satisfying $\stp(b''/C)=\stp(b/C)$ and $\fg {b''} C b$, if $\rho(a',b,b'')$ holds, then $\tp(a/C)=\tp(a'/C)$.
%\end{Conclusion}
%
%\bp  By Lemma~\ref{getisol}(2),  $\tp(bb''/C)=\tp(bb'/C)$ where $b'$ is chosen as there.  Thus, we finish by (3).
%\endproof
%

%More precisely, in order to establish $\fg a {\Cbar bb'} {BB'}$,
%one splits $BB'$ into its the $p$-weight zero part and its $p$-semiregular part.
%The independence of the first half is due to $(a,b,C)$ having minimal strength,
%and the independence of the second half is due to the fact that $w_p(a/Cbb')=0$
%from (1).
%\endproof

Continuing, we get an analogue of Lemma~3.4\cite{HrSh}.

\begin{Lemma}  \label{3.4}  There is a formula $\alpha(x)$ over $Cb$ such that, for any $a'$, $\alpha(a')$ implies:
\begin{enumerate}
\item  $w_p(a'/C)\le 2k$, $w_p(a'/Cb)\le k$, $w_p(b/Ca')\le k$; and
\item  For all $b'$, if $\stp(b'/Ca')=\stp(b/Ca')$ and $w_p(b'/Ca'b)\ge k$, then both $\rho(a',b,b')$ and $w_p(a'/Cbb')=0$ hold.
\end{enumerate}
\end{Lemma}

\bp  First, all of the inequalities involving $p$-weight are definable by  Lemma~\ref{pweightprov}.  In particular, (1) holds.
 As for (2), we first note that it holds for $a$.
Choose any $b'$ with $\stp(b'/Ca)=\stp(b/Ca)$ and $w_p(b'/Cab)\ge k$.  Since $w_p(b/Ca)=k$ by Proposition~\ref{k},
we conclude that $w_p(b'/Ca)=w_p(b/Ca)=k$. As both $\stp(b/Ca)$ and $\stp(b'/Ca)$ are $p$-semiregular, we conclude that $\fg {b}  {Ca} {b'}$.  Thus
$\rho(a,b,b')=0$ and $w_p(a/Cbb')=0$ hold by Lemma~\ref{getisol}.  As the implication  holds for every $b'$ and since `having the same strong type' is describable by a set of formulas, it follows by compactness that a $Cb$-definable  formula $\alpha(x)$ exists.
\endproof

Next, Lemma~3.5 of \cite{HrSh} becomes:

\begin{Lemma} \label{nextHrSh}  For any $a'$, if $w_p(a'/Cb)=k$ and  $\alpha(a')$ holds, then $\tp(a'/Cb)=\tp(a/Cb)$.
\end{Lemma}

\bp  Fix any such $a'$ and choose $b'$ such that $\stp(b'/Ca')=\stp(b/Ca')$ and $\fg {b'} {Ca'} b$.  
We first compute $w_p(a'b/C)=w_p(a'/bC)+w_p(b/C)=k+2k=3k$.   Since $w_p(a'/C)\le 2k$, we have $w_p(b/Ca')=k$.
Since $\stp(b'/Ca')=\stp(b/Ca')$, we have $w_p(b'/Ca')=k$ as well.  Thus, by Lemma~\ref{3.4} we conclude that $w_p(a/Cbb')=0$ and $\rho(a',b,b')$.
But then
$$w_p(bb'/C)=w_p(a/C)+w_p(b/Ca)+w_p(b'/Cab)-w_p(a'/Cbb')=2k+k+k-0=4k$$
As $\stp(b/C)$ is $p$-semiregular with $w_p(b/C)=2k$, the same holds for $\stp(b'/C)$, hence $\fg b {C} {b'}$.  
Thus, $\tp(a'/Cb)=\tp(a/Cb)$ by Lemma~\ref{getisol}(4).
\endproof

We now finish the proof of Theorem~\ref{strong}.   Recall that $\tp(e/\emptyset)$
is stationary and $p$-semiregular with $w_p(e/\emptyset)=k$.
 Since $(a,b,C)$ envelops $e$,  $e \anchor {Cb}$ and there is a 0-definable function $h$ such that
$h(a)=e$.  
%Since $\stp(a/Cb)$ is $p$-semiregular of $p$-weight $k$ and $\nfg a {Cb} e$, $w_p(a/Cbh(a))\le k-1$, so by strengthening $\alpha$, we may additionally
%assume that $w_p(a'/Cbh(a'))\le k-1$ whenever $\alpha(a')$ holds.
Let $\alpha^*(x):=\exists x'(\alpha(x')\wedge h(x')=x)$.  We claim that $\tp(e/Cb)$ is $k$-strongly regular via $\alpha^*(x)$. 
 To see this, suppose $\alpha^*(e')$ holds and $w_p(e'/Cb)=k$.
Choose $a'\in h^{-1}(e')$ with $\alpha(a')$ holding.  Then $w_p(a'/Cb)=k$ as well, so $\tp(a'/Cb)=\tp(a/Cb)$ by Lemma~\ref{nextHrSh}.  In particular, $\tp(e'/Cb)=\tp(e/Cb)$,
so   
$\tp(e/Cb)$ is $k$-strongly regular.
In light of Remark~\ref{add}, so is $\tp(e/\emptyset)$.
\endproof

%\begin{itemize}
%\item  In (1), multiply all inequalities by $k$;
%\item  Replace (2)ii by `$w_p(b'/Cab)<k$';
%\item  (2)iii remains as stated, `$\rho(a,b,b')$ holds and $w_p(a/Cbb')=0$.'
%\end{itemize}
%
%With these changes, Lemma~3.5\cite{HrSh} goes through.  Finally, the Proof of ~3.1\cite{HrSh} goes through verbatim, noting our definition of $(a,b,C)$ being `over $e$' implies
%$\fg e {\phantom{e}} {Cb}$.  In particular, $\stp(e/Cb)$ is strongly $p$-semiregular.  In light of Remark~\ref{add}, so is $\stp(e/\emptyset)$.
%\endproof

\section{Applications of Theorem~\ref{strong}}

In this brief section, we give three applications of Theorem~\ref{strong},  although they will not be used in the proof of our main results.
%
%\medskip
%\centerline{{\bf For the whole of this section we assume $T$ is (at least) stable.}}
%\medskip
%

\begin{Lemma}  \label{stray}  Suppose $T$  is superstable.  If $H$ is an infinitely definable connected group (over $A$) whose   generic 
type is strongly $p$-semiregular, then $H$ is definable over $A$.
\end{Lemma}

\bp Let $r$ be the generic type of $H$ and $\phi(x) \in r$ be a formula such that $r$ is the unique type of $p$-weight $k$ in $\phi$. By superstability, there is a definable group $H'$ with connected component $H$. The formula $\phi(x)$ contains the principal  generic, hence $H'$ is a finite union of translates of $\phi$. Each translate of $\phi$ contains a unique type of $p$-weight equal to $k$. Every other generic type of $H'$ must also be of $p$-weight $k$, so there are only finitely many generic types in $H'$. This means that $H$, its connected component,   has finite index in $H'$, hence that $H$ itself is in fact definable ($H$ is a closed subgroup in $H'$, if it has finite index it must be also open).
\endproof

Together with  Theorem~\ref{strong}, this immediately yields the following corollary.

\begin{Corollary} Let $T$ be classifiable and $p$ a non-locally modular regular type. If $H$ is an infinitely definable connected group (over $A$), with  generic $p$-semiregular, then $H$ is definable over $A$. \end{Corollary}

\medskip

The second application connects binding group constructions with isolation of $p$-semiregular types, in the spirit of \cite{HrShOT} or more recently \cite{MoosaSanchez}.
First,  we record an easy remark.

\begin{Remark} {\em In general, if there is a $B$-definable  group $G$ with a $B$-definable transitive action on a complete type $q$ over $B$, then $q$ is isolated over $B$: Let $e$ realize $q$, consider the formula $\phi(x) \in L(Be)$ asserting $\exists g \in G (x = g.e)$. Then $\phi(x) $ holds if and only if $x \models q$. So $q$ is isolated by a formula over $Be$. On the other hand, $q$ is $B$-invariant, so in fact, $q$ is isolated by a formula over $B$. }
\end{Remark}

Let us now recall some definitions and notation.

\begin{Definition} {\em Let $B = acl (B)$ and let $p,q\in S(B)$ with $p$ regular.  
We say that $q$ is {\em $p$-internal} if there is $D\supseteq B$ such that for some $a\models q|D$, 
$a\in\dcl(D\cup p(\C))$.
}
%and $q$ be a (partial)  type over $B$ and $r$ a complete type over $B$. We say that $r$ is {\em $q$-internal} over $B$ if there is $D\supset B$ such that for every $e \models r|D$,  $e\in dcl (q(\C)$)}.
\end{Definition}

The following existence result is Corollary 7.4.6 of \cite{Pillay}.

\begin{Proposition} \label{internal}  Let $T$ be stable.  Suppose $B$ is algebraically closed and $p,q\in S(B)$ are non-orthogonal with $p$ regular.
Then for any $a$ realizing $q$,  there is $a'\in\dcl(Ba)\setminus B$ such that $\tp (a'/B)$ is $p$-internal.
\end{Proposition}

\begin{Lemma}\label{indepforsr} Let $T$ be stable.  Suppose $B=\acl(B)\supseteq M$, $p\in S(M)$ regular and $q\in S(B)$ is $p$-semiregular.  If $q\perp^a M$ then
   $\fg a B {p(\C)}$ for any $a$ realizing $q$.
\end{Lemma}

\bp Let $a$ realize $q$ and let $E\subseteq p(\C)$ be finite.  Let $E_1\subseteq E$ be a maximal $M$-independent subset of $E$, and let
$E_2$ be a maximal subset of $E_1$ satisfying $\fg {E_2} M B$.
Note that  $E_2$ realizes a Morley sequence $(p|B)^{(k)}$ of the non-forking extension of $p$ to $B$.

\medskip
\noindent
{\bf Claim.}  $w_p(E/E_2B)=0$.

\medskip
\bp  First, choose any $e\in E_1\setminus E_2$.  Because maximality implies that $\nfg {eE_2} M B$, $\tp(e/BE_2)$ is a forking extension of a regular type, hence $w_p(e/E_2B)=0$.
It follows that $w_p(E_1/E_2B)=0$ as well. 
Now choose any $e\in E\setminus E_1$.  By the maximality of $E_1$, $\nfg e M {E_1}$ so $w_p(e/E_1)=0$.  Thus, $w_p(e/E_1B)=0$.  As this holds for every $e\in E\setminus E_1$,
$w_p(E/E_1B)=0$.  Combining the two arguments yields $w_p(E/E_2B)=0$.
\endproof

Because $\tp(E_2/B)$ does not fork over $M$ and $q\perp^a M$, we have that $\fg a B {E_2}$.  As $q$ is $p$-semiregular, that $\fg a B E$ follows immediately from the Claim.\endproof

\begin{Theorem}\label{bindinggroup}{\em (Binding group)}  {\em(\cite{Pillay} 7.4.8 or \cite{Poizat} 2.2.20)}   
Suppose $T$ is stable.  Let $B = acl(B)$ and let $p,q\in S(B)$  with $p$ regular and $q$ $p$-internal. Let $G := Aut (q(\C) / B \cup p(\C))$. Then $G$ and its action are infinitely definable in the following sense: there is $G_1$ an infinitely definable group over $B$ and a $B$-definable action of $G_1$ on $q(\C)$ such that, as permutation groups of $q(\C)$ over $B\cup p(\C)$, $G$ and $G_1$ are isomorphic.
Additionally, if $\fg a B {p(\C)}$ for all $a$ realizing $q$, then $G$ and hence also $G_1$ act transitively on $q(\C)$. 
\end{Theorem}

Note that by stability, the group $G$ is   also the group of restrictions to $q(\C)$ of the automorphisms of $\C$ fixing $B\cup p(\C)$ pointwise. $G$ is infinite if and only if $q $ is not algebraic over $B\cup p(\C)$.
%We can always suppose that the action is faithful (that is, that the subgroup of $G_1$ which fixes all of $q(\C) $ pointwise is trivial). 
The next proposition tells us that we get a transitive action from the connected component of $G_1$.

\begin{Proposition} If an infinitely definable group $G$ is defined over $B= acl(B)$ and there is a $B$-definable, transitive action of $G$ 
on a complete type  $q\in S(B)$,
then the restriction of the action to the connected component $G^0$ of $G$ is also transitive.  
\end{Proposition}

\bp For notational simplicity suppose $B =\emptyset$.

%\noindent As $G^0$ is a subgroup of $G$, everything is clear except that $G^0$ acts transitively on $Q$.  We obtain this via two claims and a brief argument.

\medskip
\noindent
{\bf Claim 1.}  There is some pair $(e,e^*)$ realizing $q\otimes q$ and some $h\in G^0$ such that $h.e=e^*$.

\bp  Choose a set of representatives $R\subseteq G(\C)$ such that every $g\in G$ can be written as $ch$ for some $c\in R$ and $h\in G^0$, i.e., $R$ contains an element of every $G^0$-coset of $G$.  As the index $[G:G^0]$ is bounded, we may choose $R$ of bounded size ($2^{\aleph_0}$ if the language is countable).
Choose any $e$ realizing $q|R$ and any $e'$ realizing $q|Re$.

By the transitivity of the action, as both $e,e'$ realize $q$, choose $g\in G$ such that $g.e=e'$.  By choice of $R$, choose $c\in R$ and $h\in G^0$ such that
$g=ch$ and put $e':=h.e$.  We need to show that $e^*\anchor e$.  For this, note that
$$e^*=(c^{-1}g).e=c^{-1}.e'\in\dcl(Re')$$
But $e'\anchor Re$ and $e\anchor R$ imply $e'R\anchor e$, hence $e^*\anchor e$ as required.
\endproof
%
%Now, as $c\in R$, both $g$ and $h$ are equi-definable over $R$.  Thus, $g.e$ and $h.e$ are equi-definable over $Re$.  As $g.e=e'$
%and $\fg {e'} {\emptyset} {Re}$, we conclude that $\fg {h.e} {\emptyset} {Re}$.
%Thus, $e^*:=h.e$ satisfies the Claim.
%\endproof

\medskip
\noindent
{\bf Claim 2.}  For every $(e,e')$ realizing $q\otimes q$ there is $h\in G^0$ such that $h.e=e'$.

\bp  Homogeneity of $\C/B$.  Fix $(e,e^*)$ and $h$ as in Claim 1, and let $(e_1,e_2)$ be any other realization of $q\otimes q$.
Choose an automorphism $\sigma$ of $\C$, fixing $B$ pointwise with $\sigma(e)=e_1$ and $\sigma(e^*)=e_2$.
As $G^0$ and the action are $B$-definable, $\sigma(h)\in G^0$ and $\sigma(h).e_1=e_2$.
\endproof

To complete the proof of the Proposition, choose any $e,f\in q(\C)$.  Choose $e^*$ realizing $q|\{e,f\}$.  By Claim 2, choose
$h_1\in G^0$ such that $h_1.e=e^*$ and choose $h_2\in G^0$ such that $h_2.e^*=f$.  Then $h_2h_1\in G^0$ and
$(h_2h_1).e=f$, so $G^0$ acts transitively on $q(\C)$.
\endproof

The following results give a sufficient condition for a non-locally modular type to be isolated. These results will be used as part  of the  forthcoming work of the authors on the analysis of weight one models in classifiable theories.
The following definition appears already in \cite{HrShOT}.

\begin{Definition} {\em   Suppose $T$ is superstable.  We say that $q \in S(B)$ is {\em $c$-isolated} (following Hrushovski-Shelah in \cite{HrShOT}) if
  there is a formula $\theta(x)\in q$ such that $R^\infty (\theta(x)) = R^\infty (q) = \alpha$ and furthermore, for any $r \in S(B)$, $\theta(x)\in r$ implies
  $R^\infty (r) = \alpha$. 
  }
  \end{Definition}

\begin{Proposition}\label{groupisol} Suppose $T$ is superstable.  Let $q \in S(B) $ be $p$-strongly semiregular and $c$-isolated via the formula $\phi(x)$ and let $G_1$ be an infinitely definable group over $B$, with a 
$B$-definable, transitive action of $G_1$ on $q(\C)$. Then $q$ is isolated. 
\end{Proposition}

\bp  By the usual construction (e.g., Lemma~1.6.19 in \cite{Pillay}) find a $B$-definable supergroup $H\supseteq G_1$, a $B$-definable set $X$ containing $q$, and a $B$-definable,  transitive action of $H$ on $X$ which extends the action of $G_1$ on $q(\C)$, and such that   $X \subset \phi(\C) $.

  So $X$ has the property that every type over $B$ in $X$ has the same $R^\infty$ rank as $q$, say $\alpha$, that every type in $X$ is $p$-simple, of $p$-weight at most equal to $k $, the $p$-weight of $q$, and that $q$ is the unique type in $X$ of $p$-weight exactly  $k$. We show that $X$ isolates $q$.
  
  \medskip\noindent{\bf Claim 1.}  Let $e $ realize $q$, and let $h \in H$ be independent from $e$ over $B$, then $h.e$ also realizes $q$. 
  
%  \bp  Then $e$ and $h,h^{-1}$ are also independent over $B$. If $d = h.e$, then $d$ and $e$ are interdefinable over $B,h,h^{-1}$. As  $X \subset \phi(\C) $, $t(d/B)$ is $p$-simple of $p$-weight at most $k$. As $e$ is independent from $h,h^{-1}$ over $B$, it remains of $p$-weight $k$. By interdefinablity, $d$ and $e$ must have same $p$-weight over $B,h,h^{-1}$, that is $k$. It follows that $t(d/B)$ must already have $p$-weight $k$ over $B$, and hence must be equal to $q$.
%  \endproof
\bp  Letting $d=h.e$, $d$ and $e$ are inter-definable over $Bh$.  Since $e$ realizes $q|Bh$, $w_p(e/Bh)=k$, hence $w_p(d/B)\ge k$.
By our assumptions on $X$, $w_p(d/B)=k$, so $d$ realizes $q$.
 \endproof

 \medskip\noindent{\bf Claim 2.} 
 Let $h$ be any element of $H$, not necessarily independent from $e$, and let $d = h.e$. Then $d$ realizes $q$.
 
 \bp
  Let $g$ be a generic in $H$, independent from $h,e$ (and hence from $h,e,d$) over $B$. It follows that  $g^{-1}h $ and $e$ are independent over $B$ : $\fg {g^{-1}} {Bh} e$ by choice of $g$, hence $\fg {(g^{-1} h)} {Bh} e$. As $g^{-1}$ is generic and independent from $h$, $g^{-1} h$ is also  independent from $h$ (see e.g., 1.6.9(iv) of \cite{Pillay}), so it follows that $g^{-1}h $ and $e$  are independent over $B$.
  Hence by Claim 1,  $f := (g^{-1} h).e$ realizes $q$.
  
We next show $g $ and $f$ are independent over $B$: Since  $f = g^{-1}. d$, $d$ and $f$ are inter-definable over $Bg$, hence they have same $\infty$-rank over $Bg$.  
By $c$-isolation, $R^\infty(d/B)=\alpha$.  By our choice of $g$, $\fg g {B} d$, hence $\alpha=R^\infty(d/B) = R^\infty(d/Bg) = R^\infty(f/Bg)$.  As $f\in X$, $R^\infty(f/B)\le\alpha$, so 
$R^\infty(f/Bg)=R^\infty(f/B)$, implying $f$ and $g$ are independent over $B$. 

Thus, by Claim 1, $g.f$ must realize $q$, but $g.f = d$.
  \endproof

%To finish, we show $X\subseteq q(\C)$.  Choose any $d\in X$ and $e\in q(\C)\subseteq X$.  As the action of $H$ is transitive, choose $h\in H$ with $h.e=d$.
%By Claim 2, $d\in q(\C)$.
  As the action of $H$ on $X$ is transitive,  
 % every $d\in X$ can be written as $d=h.e$ for some $h\in H$, $e\in q(\C)$.
  Claim 2 implies that any $d\in X$ must realize $q$.   That is, $X$ isolates $q$.
  \endproof

  \begin{Corollary} \label{point} Let $T$ be classifiable with $M\subseteq_{na}\C$ and $B\supseteq M$ algebraically closed and suppose
   $p$ is a non-locally modular regular type with $p\not\perp M$.  If $q=\tp(a/B)$ is $p$-semiregular, $c$-isolated, but $q\perp^a M$,
   then there is $a'\in\dcl(Ba)\setminus B$ such that $q':=\tp(a'/B)$ is isolated.
   \end{Corollary}

  \bp  Since $M\subseteq_{na}\C$  and $p\not\perp M$ is regular,   there is a regular $p_0\in S(M)$ non-orthogonal to $p$ by Fact~\ref{qregularpair}.
  To ease notation, we may assume that our $p$ is the non-forking extension of $p_0$ to $S(B)$.
 % Since $q$ is $p$-semiregular and $q\not\perp M$, $p\not\perp M$.  Thus, as $M\subseteq_{na} \C$, by XYZ there is a regular $p_0\in S(M)$ non-orthogonal to $p$.
  %Let $p_1$ denote the non-forking extension of $p_0$ to $S(B)$ and choose any realization $a$ of $q$.  
   By Proposition~\ref{internal} choose $a'\in\dcl(Ba)\setminus B$ such that $q':=\tp(a'/B)$ is $p$-internal.   
   It is easily checked that $q'$ remains $c$-isolated and $p$-semiregular and that $q'\perp^a M$.  
   Since  $T$ is classifiable and $p$ is non-locally modular, $q'$ is strongly $p$-semiregular by Theorem~\ref{strong}.
   
  As $q\perp^a M$, by Lemma~\ref{indepforsr} we have
  $\fg {a'} B {p_0(\C)}$, hence $\fg {a'} B {p(\C)}$ as well.  Thus, by Lemma~\ref{bindinggroup}, there is an infinitely $B$-definable group $G_1$ with a $B$-definable, transitive action
  on $q'(\C)$.  Hence $q'$ is isolated  by Proposition~\ref{groupisol}.
  \endproof

%  $B=\acl(B)\supseteq M$.  Suppose $q,p\in S(B)$ with $p$ regular and non-locally modular and $q$ $p$-internal, $p$-semiregular, and $c$-isolated.
%  $p$ regular non-locally modular. If $q\in S(B)$ ($B = acl(B)$) is $p$-semiregular, $c$-isolated and $q\perp^a M$, then $q$ is isolated.\end{Corollary}
%  \bp First by Theorem~\ref{strong}, $q$ is $p$-strongly semiregular. By classifiability, $q \not\perp M$ (Fact~\ref{NDOP}); by  Proposition \ref{internal} and the binding group Theorem (Proposition~\ref{bindinggroup}), there is an infinitely definable group $G$, defined over $B$ which acts transitively on $q(\C)$. It now follows by Proposition~\ref{groupisol} that $q$ is isolated. 
%  \endproof

 Note that in  Corollary~\ref{point} we used two conditions on the type $q'$ to prove isolation: the strong $p$-semiregularity condition goes up to non forking extensions, but the $c$-isolation does not necessarily.

We finish this section with a third   way of  proving isolation, without group actions, in the case of a strongly regular type. It is not clear if this 
method could be generalized to the case of strong $p$-semiregularity.

\begin{Proposition}\label{Heresy}  Suppose that $T$ superstable, $M\subseteq_{na}\C$, $B\supseteq M$ is algebraically closed, and $q\in S(B)$ is
 strongly regular, depth zero, and that $q \perp^a M$. Then either $q \perp M$ or $q$ is isolated. 
\end{Proposition}

\bp  If $q$ is trivial, $q\perp^a M$ implies that $q\perp M$ (see, e.g., X 7.3(6) of \cite{Shc}). So we can suppose that $q$ is non-trivial.
Henceforth, assume that $q$ is  non-orthogonal to $M$. By Fact~\ref{qregularpair} choose a regular $p\in S(M)$ with
$q\not\perp p$.   As $q\perp^a M$, in particular, over $B$, $q$ is almost orthogonal to 
$p^{(\omega)}$.
As $q$, and hence $p$, is non-trivial of depth 0, use  Fact~\ref{pweightprov} to 
choose $\phi\in q$ that is $p$-simple, such that $p$-weight is defined and continuous in $\phi$.  As $q$ is strongly regular,
by strengthening $\phi$
we may additionally assume $q$ is the only type over $B$ containing $\phi$ of positive $p$-weight.

Now  choose $n$ least such that there are $\abar=(a_1,\dots, a_n)$ with each $a_i$ realizing $\phi$ and $tp(\abar/B)$ is not almost orthogonal to $p^{(\omega)}$.
We know that such a finite $n$ exists, since $q$ non-orthogonal to $p$ implies that some  $q^{(\ell)}$ is not almost orthogonal to $p^{(\omega)}$ (see e.g., 4.3.1(iii) of \cite{Pillay}). Remark:  Here, however, we are minimizing $n$
without assuming $\abar$ is $B$-independent.  

Note that $n\ge 2$.  Indeed,  if $a_1$ realizes $q$, then $a_1/B$ is almost  orthogonal to $p^{(\omega)}$ by assumption.  
On the other hand, if $a_1$ realizes $\phi$ but not $q$, then $w_p(a_1/B)=0$, so $a_1$ cannot fork  over $B$ with any any independent set of realizations of $p$.

Once $n$ is fixed, choose $k$ such that $\abar/B$ is not almost orthogonal to $p^{(k)}$.  To save writing, let $n=m+1$ and $r(\ybar):=(p|B)^{(k)}$.
Choose a specific realization $\cbar$ of $r$ such that $\nfg {\abar} B {\cbar}$, and choose an $L(B)$-formula $\theta(\xbar,\ybar)\in\tp(\abar\cbar/B)$ witnessing the forking.
Let
$$\gamma(x_0):=\phi(x_0)\wedge d_r\ybar \left[\exists x_1\dots\exists x_m(\bigwedge \phi(x_i)\wedge\theta(x_0,x_1,\dots,x_m,\ybar))\right]$$
As $B$ is algebraically closed, $\gamma$ is over $B$.  We argue that $\gamma$ isolates $q$.

To see this, choose any $b_0$ realizing $\gamma$.  Choose $\dbar$ realizing $r|Bb_0$, and choose witnesses $b_1,\dots,b_m$.
Thus, the $n$ elements $b_0,\dots,b_m$ each realize $\phi$ and $\theta(\bbar,\dbar)$ holds.   We argue that in fact, every $b_i$ realizes $q$.
Let $I=\{i\le m: b_i$ realizes $q\}$.  Let $\bbar_I$ be the subsequence of $\bbar$ induced by $I$.
By way of contradiction, assume  $|I|<m+1=n$.  By the minimality of $n$, we must have $\fg {\bbar_I} B {\dbar}$.  But also, by our choice of $\phi$, if $i\not\in I$, then
$w_p(b_i/B)=0$.  Thus, $w_p(\bbar/B\bbar_I)=0$.  But $\dbar$ is a Morley sequence in $p$, hence $\stp(\dbar/B\bbar_I)$ is $p$-semiregular.
Thus, $\fg {\bbar}  {B\bbar_I} {\dbar}$.  It follows by transitivity that $\fg {\bbar} B {\dbar}$, which is contradicted by $\theta(\bbar,\dbar)$.
\endproof

\begin{Corollary} \label{Heresyclassif} Let $T$ be classifiable, $M\na \C$, $B\supseteq M$ algebraically closed,  and let $q\in S(B)$ be regular, but  non locally modular. 
If  $q\not\perp M$ but $q \perp^a M$, then $q$ is isolated. 
\end{Corollary}
\bp By the classifiability of $T$, $q$ must also  be strongly regular of depth zero, so Proposition~\ref{Heresy}  applies.\endproof

\section{Constructible, minimal models over realizations of non-locally modular types}

Recall our global assumption that all theories $T$ are (at least) stable.
The basic definitions and facts about locally modular regular types can be found in the Appendix (Section \ref{appendix-locallymodular}).

\subsection{Witnesses to non-modularity}

We describe a minimal counterexample to non-modularity in our context.  

\begin{Definition} {\em  Let $M$ be any model (not necessarily an a-model) and let $p\in S(M)$ be regular.
A quadruple $(a,b,c,d)$ is {\em 4-dependent} if
\begin{enumerate}
\item  Any three elements realize $p^{(3)}$, but
\item  $w_p(abcd/M)=3$.
\end{enumerate}
A {\em witness to non-modularity for $p$ over $M$} is a set of parallel lines, i.e.,  some 4-dependent quadruple $(a,b,c,d)$ such that
$\clp(Mab)\cap\clp(Mcd)=\clp(M)$.
}
\end{Definition}

The following proposition follows from Section \ref{appendix-locallymodular} of the appendix. 
\begin{Proposition}\label{amodelok} Over any a-model $M$,  if $p\in S(M) $ is a regular type, then $p$ is non-locally modular if and only if there is a witness to non-modularity for $p$ over $M$.
\end{Proposition}

\begin{Lemma} \label{exists}
 Suppose $M$ is any countable model, and $p\in S(M)$ is regular, but not locally modular.   There is an $\epsilon$-finite set $E$ ($E$ is contained in the algebraic closure of some finite set) such that for any
countable model $M'$ containing $M\cup E$, 
there is a witness $(a,b,c,d)$ to non-modularity over  $M'$.
\end{Lemma}  

\bp  Let $M^*\supseteq M$ be any a-model, and let $q$ denote the non-forking extension of $p$ to $M^*$.  
As $q$ is not locally modular, it follows from Proposition~\ref{amodelok} that there is a witness $(a,b,c,d)$ to non-modularity over  $M^*$.  Choose a countable, $\epsilon$-finite $E$ over which $\tp(abcd/M^*)$ is based and stationary.  To see that $E$ suffices, 
choose  any countable $M'$ containing $M\cup E$.  
As $M^*$ is sufficiently saturated, we may assume that $M'\preceq M^*$.
We argue that $(a,b,c,d)$ is
a witness to non-modularity over  $M'$.  To see this, note that $\fg {abcd} {M'} {M^*}$.  Thus, $(a,b,c,d)$ is 4-dependent
with respect to the non-forking extension of $p$ to $M'$.
For the final clause, let $C:=\dcl(M'abcd)\cap\clp(M')$.
Note that $\fg {abcd} C {M^*}$ and
$\stp(abcd/C)$ is $p$-semiregular by Criterion~\ref{Criterion}. Thus,  $\clp(M^*ab)\cap\clp(M^*cd)=\clp(M^*)$ implies that
$\clp(Cab)\cap\clp(Ccd)=\clp(C)$ by Proposition~\ref{nonforking}.  However, $\clp(C)=\clp(M')$, $\clp(Cab)=\clp(M'ab)$,
and $\clp(Ccd)=\clp(M'cd)$, so we finish.
\endproof

\subsection{A definable witness to non-modularity}

In this section we assume throughout that $T$ is classifiable.    Under this assumption, we obtain a definable   witness to non-modularity
over an arbitrary model $M$ on which a non-locally modular type is based.  The starting point is the following, which melds Theorem~1b of \cite{HrSh} with
Theorems~\ref{pweightprov} and \ref{NDOP} from the appendix.

%By Lemma~\ref{exists}, if $p\in S(B)$ is any non-locally modular
%stationary, regular type  over a countable set $B$, then there is a countable $M$ containing $B$ with a witness to non-modularity $(a,b,c,d)$
%over $M$ (relative to the non-forking extension of $p$ to $M$).  We first explore how definable  such a witness is. For this, note that by 
%Theorems~\ref{pweightprov} and \ref{NDOP} 
%we have:
%%some theorems of Hrushovski and Shelah in \cite{HrSh}.

\begin{Theorem} \label{olddefinability}
If $T$ is classifiable, $B$ is algebraically closed, and  $p\in S(B)$ is a regular, non-locally modular type, then there is a formula 
$\theta\in p$ such that (recall Definition~\ref{defpsimple}):
\begin{enumerate}  
\item $\theta$ is $p$-simple of $p$-weight one;

\item  $p$-weight is defined and continuous inside $\theta$; and

\item $p$ is strongly regular via the formula $\theta$, i.e., for every $C\supseteq B$ and every $e\in\theta(\C)$, if $w_p(e/C)>0$, then $\tp(e/C)$ is the non-forking extension of $p$ to $C$.
\end{enumerate}
\end{Theorem}

For the next few Lemmas, fix a model $M$ and non-locally modular $p\in S(M)$.  Also fix $\theta(x,\mbar)\in p$ as in Theorem~\ref{olddefinability} with $\theta(x,\ybar)$ 0-definable.

\begin{Lemma}  \label{4dep} Suppose $p\in S(M)$ is regular and $(a,b,c,d)$ is 4-dependent over $M$.  
Then there is some $h\subseteq M$ and an $\mbar$-definable formula $R(x,y,z,w,u)\in \tp(abcdh)$ such that for any $(a',b',c',d')$ with some triple realizing $p^{(3)}$ and
for any $h'\subseteq M$, if $R(a',b',c',d',h')$ holds, then $(a',b',c',d')$ is 4-dependent over $M$.
\end{Lemma}

\bp  First, as $w_p(a/bcdM)=0$ and $\theta(a)$ holds, choose  $h_1\subseteq M$ and $\alpha_1(x,y,z,w,u)\in\tp(abcdh_1)$ such that for any 
$a',b',c',d'$
and any $h_1'\subseteq M$, if $\alpha_1(a',b',c',d',h_1')$ holds, then  $\theta(a')$ holds and $w_p(a'/b'c'd'M)=0$.
Similarly, choose $h_2,h_3,h_4\subseteq M$ and $\alpha_2,\alpha_3,\alpha_4$ such that e.g., $\alpha_2(a',b',c',d',h_2')$ implies $\theta(b')$ and 
$w_p(b'/a'c'd'M)=0$.  Put $h:=h_1h_2h_3h_4$ and put $R(x,y,z,w,u):=\bigwedge_{i=1}^4 \alpha_i$.
   
We argue that $R$ is as claimed.  By symmetry, choose $b'c'd'$ realizing $p^{(3)}$ and $h'\subseteq M$.  We verify that for every $a'\in \C$, if
$R(a',b',c',d',h')$ holds, then $(a',b',c',d')$ is 4-dependent over $M$.  To see this, since $\alpha_1(a',b',c',d',h')$ holds, we have $\theta(a')$  and $w_p(a'/b'c'd'M)=0$.
Because of $\alpha_2$ holding, $w_p(b'/a'c'd'M)=0$.  As $b'c'd'$ realizes $p^{(3)}$, we must have $w_p(a'/M)>0$.  Thus, by strong regularity,
$\tp(a'/M)=p$.  To verify 4-dependence over $M$, it remains to show that each of $a'b'c'$, $a'b'd'$ and $a'c'd'$ realize $p^{(3)}$.  We check this for $a'b'c'$, with the other two following by symmetry.  As $\alpha_4(a',b',c',d',h')$ holds, $w_p(d'/a'b'c'd'h')=0$.  But, as $b'c'd'$ realizes $p^{(3)}$ this implies $\nfg {d'} {Mb'c'} {a'}$.  As $\tp(d'/Mb'c')$ is regular, this implies
$\fg {a'} M {b'c'}$, so $a'b'c'$ realizes $p^{(3)}$.
\endproof

%\begin{Lemma} \label{4dd}
%Suppose that $p\in S(M)$ is regular and $(a,b,c,d)$ is a  4-dependent sequence of realizations of $p$.  Then there is a symmetric formula
%$S(x_1,\dots,x_4)\in\tp(abcd/M)$ such that for any $(b',c',d')$ realizing $p^{(3)}$ and for any element $a'$ realizing $S(x_1,b',c',d')$,
%we have $(a',b',c',d')$ is 4-dependent.
%\end{Lemma}
%
%\bp  For a given ordering of the variables, let $S_0$ be the conjunction of $\bigwedge_{i=1}^4\theta(x_i)$ with formulas in $\tp(a,b,c,d)$ demonstrating
%that each variable has $p$-weight zero over the other three.  Then, for any three coefficients, which we write as $(x_2,x_3,x_4)$
%for definiteness, if $(b',c',d')$ realizes $p^{(3)}$ and $S_0(a',b',c',d')$ holds, then $\tp(d'/Mb'c'a')$ forks over $Mb'c'$.  As $\tp(d'/Mb'c')$
%is a non-forking extension of $p$, this implies $w_p(a'/M)>0$.  Combined with $\theta(a')$, this implies that $\tp(a'/M)=p$.  Easy $p$-weight computations
%show that $(a',b',c',d')$ is 4-dependent.  To find a symmetric $S$, take the disjunction of the 24 $S_0$'s, with respect to each ordering of
%$(x_1,\dots,x_4)$.
%\endproof

Next, among all 4-dependent quadruples $(a,b,c,d)$, we want to distinguish those that are witnesses to non-modularity over $M$.
This is the content of the next two Lemmas.

\begin{Lemma}  \label{firstL}  Suppose $p\in S(M)$ is a regular, non-locally modular type and $(a,b,c,d)$ is a witness to non-modularity over $M$.
Then for every $e$ realizing $p|Mabcd$ and for every $f\in\theta(\C)$ that satisfies $b\in\clp(Maef)$, we have $(c,d,e,f)$
realizes $p^{(4)}$ and $w_p(ab/Mcdef)=0$.
\end{Lemma}

\bp  First, since $\theta(f)$ holds, $\tp(f/M)$ is $p$-simple.  As $b$ realizes $p$, $\fg b M {ae}$, but $b\in\clp(Maef)$, we must have
$w_p(f/M)>0$, hence $\tp(f/M)=p$.  Thus, all six elements $a,b,c,d,e,f$ realize $p$.  Next, note that $w_p(abcdef/M)=4$, since
$w_p(abcd/M)=3$ (by 4-dependence), $\fg e M {abcd}$, and $w_p(f/Mabcde)=0$ by exchange.
Next, we show that $w_p(cdef/M)=4$, i.e., that $(c,d,e,f)$ realizes $p^{(4)}$.
By way of contradiction, assume that this were not the case, i.e, that $w_p(cdef/M)\le 3$.  
We compute the following $p$-weights:
\begin{itemize}
\item  $w_p(ef/Mabcd)=1$ [it is $\ge 1$ because of $e$, but $<2$ since $f\in\cl_p(Mabe)$].
\item  $w_p(ef/Mcd)=1$ [it is $\ge 1$ from the former line, but if it was $=2$, then we would have $w_p(cdef/M)=4$].
\item  $w_p(ef/Mab)=1$ [it is $>0$ because of $e$, but $<2$ because $f\in\cl_p(Mabe)$].
\end{itemize}
As both $\stp(ef/\clp(Mab))$ and $\stp(ef/\clp(Mcd))$ are $p$-semiregular, 
$$\fg {ef} {\clp(Mab)} {abcd}\qquad \hbox{and}\qquad \fg {ef} {\clp(Mcd)} {abcd}$$
Hence, 
$$Cb(ef/Mabcd)\subseteq\clp(Mab)\cap\clp(Mcd)=\clp(M)$$
with the last equality holding since $(a,b,c,d)$ is a witness to non-modularity.  This would imply $\fg {ef} {\clp(M)} {abcd}$,
which contradicts $f\in\clp(Mabe)$.

Finally, since $w_p(abcdef/M)=w_p(cdef/M)=4$, it follows immediately that $w_p(ab/Mcdef)=0$.
\endproof

By contrast:

\begin{Lemma}  \label{secondL} Suppose $(a,b,c,d)$ is 4-dependent over $M$, but $\clp(Mab)\cap\clp(Mcd)\neq \clp(M)$.
Then there are $ef$ such that $\stp(ef/Mab)=\stp(cd/Mab)$ with $e$ realizing $p|Mabcd$,  but $w_p(cdef/M)=3$.
\end{Lemma}

\bp  Choose $g\in (\clp(Mab)\cap\clp(Mcd))\setminus\clp(M)$.  
Since $g\in\clp(Mab)$, $\tp(g/M)$ is $p$-simple.  As $g\not\in\clp(M)$,  $w_p(g/M)>0$, but since $c\not\in\clp(Mab)$, we cannot have $w_p(g/M)\ge 2$.
Then $w_p(g/M)=1$ and $w_p(cd/Mg)=1$.  
Now choose any $ef$ such that $\stp(ef/Mabg)=\stp(cd/Mabg)$ with $\fg {ef}  {Mabg} {cd}$.  Then $e$ realizes $p|Mabcd$ and $w_p(ef/Mg)=1$.
So $$w_p(cdef/M)\le w_p(cdefg/M)\le w_p(cd/Mg)+w_p(ef/Mg)+w_p(g/M)=3$$
\endproof

%For the Corollary that follows, note that if $(a,b,c,d)$ is 4-dependent, then there is an $M$-definable formula
%$\beta(x_2,a,c,d)\in\tp(b/Macd)$ that implies both $\theta(x_2)$ and $w_p(x_2/Macd)=0$.

\begin{Lemma}  \label{defwitness}
Suppose $p\in S(M)$ is regular and that
$(a,b,c,d)$ witnesses non-modularity over $M$.  Fix an $L(M)$-formula  $\theta(x,\mbar)\in p$ from Theorem~\ref{olddefinability}.
Then there is some $h\subseteq M$ and an $\mbar$-definable  $R^*(x,y,z,w,u)\in \tp(a,b,c,d,h)$ such that, for any
$h'\subseteq M$,  any $b'c'd'$ realizing $p^{(3)}$, and any $a'\in\C$, if $R^*(a',b',c',d',h')$ holds, then $(a',b',c',d')$ witnesses the non-modularity of $p$ over $M$.
\end{Lemma}

\bp
As $(a,b,c,d)$ is 4-dependent over $M$ choose $h\subseteq M$ and $R(x,y,z,w,u)$ as in
Lemma~\ref{4dep}.  Choose any  $e\models p|Mabcd$.  Then, for    every $f$,
if $R(a,b,e,f,h)$ holds, then since $w_p(b/aefM)=0$ 
we have    $a\in \clp(Mcdef)$ by Lemma~\ref{firstL}.  Thus, by compactness (and since $p$-weight 0 formulas are closed under finite disjunctions)
there is an $L$-formula $\gamma(x,y,z,w,v_1,v_2,u)$ such that 
\begin{itemize}
\item  $\gamma(a,b,c,d,e,f,h)$  holds for any $e\models p|Mabcd$ and any $f$ for which $R(a,b,e,f,h)$ holds; and
\item  For all $(a',b,c',d',e',f',h')$,  $\gamma(a',b',c',d',e',f',h')$ holding implies $w_p(a'/c'd'e'f'h')=0$.
\end{itemize}
Put
$R^*(x,y,z,w,u):=$
$$R(x,y,z,w,u)\wedge d_pv_1 \forall v_2 [R(x,y,v_1,v_2,u)\rightarrow \gamma(x,y,z,w,v_1,v_2,u)]$$
From above, $R^*\in\tp(a,b,c,d,h)$.  Now choose any $h'\subseteq M$ and  $(a',b',c',d')$ such that $R^*(a',b',c',d',h')$ holds.
By Lemma~\ref{4dep} we know that $(a',b',c',d')$ is 4-dependent over $M$.  
By way of contradiction, suppose that $(a',b',c',d')$ is not a witness to the non-modularity of $M$.  Choose any $e'\models p|Ma'b'c'd'$.  By Lemma~\ref{secondL},
there is $f'$ such that $\stp(e'f'/Ma'b')=\stp(c'd'/Ma'b')$ but $w_p(c'd'e'f'/M)=3$.  Since $R(a',b',c',d',h')$ holds, we have $R(a',b',e',f',h')$ holding as well.    
Thus, $\gamma(a',b',c',d',e',f',h')$ holds, hence $w_p(a'/c'd'e'f'M)=0$.   By $R(a',b',e',f',h')$ again, we have $w_p(b'/a'e'f'M)=0$, so $w_p(a'b'c'd'e'f'/M)=3$.
But this is impossible since $b'c'd'$ realizes $p^{(3)}$ and $e'\models p|Mb'c'd'$.
\endproof

\begin{Proposition}  \label{overM} If $T$ is classifiable, then for any model $M$ and any $p\in S(M)$ that is regular but not locally modular,
there is a witness $(a,b,c,d)$ to the non-modularity of $p$ over $M$.  Thus, there is an $L(M)$-formula $S(x,y,z,w)\in\tp(abcd/M)$ such that
for any $b'c'd'$ realizing $p^{(3)}$, $S(x,b',c',d')$ is consistent and $(a',b',c',d')$ is a witness to non-modularity of $p$ over $M$ for every realization $a'$
of $S(x,b',c',d')$.
\end{Proposition}

\bp  It suffices to show that a witness to the non-modularity of $p$ exists, since
the second sentence follows from this via Lemma~\ref{defwitness}.
Choose $\theta(x,\mbar)\in p$ as in Theorem~\ref{olddefinability} and choose any $bcd$ realizing $p^{(3)}$.  
Choose any a-model $M^*\succeq M$ independent from $bcd$ over $M$.  By Proposition~\ref{amodelok}
there is a witness $(a_0,b_0,c_0,d_0)$ to the non-modularity of $p|M^*$ over  $M^*$.  As $\tp(bcd/M^*)=\tp(b_0c_0d_0/M^*)$, there is
some $a$ such that $(a,b,c,d)$ is a witness to the non-modularity of $p|M^*$ over  $M^*$.  
Choose an $\mbar$-definable $R^*(x,y,z,w,h)$ with $\hbar\subseteq M^*$ from Lemma~\ref{defwitness}.

So $\exists x R^*(x,b,c,d,h)$ holds.  
As  $\fg {bcd} M h$, finite satisfiability gives some $h'\subseteq M$ for which $\exists x R^*(x,b,c,d,h')$.  
Choose any $a'\in R^*(\C,b,c,d,h')$ with $\fg {a'} {Mbcd} {M^*}$.
By Lemma~\ref{defwitness} $(a',b,c,d)$ witnesses the non-modularity of $p|M^*$ over $M^*$, but also $\fg {abcd} M {M^*}$.
It follows by non-forking calculus that $(a',b,c,d)$ is 4-dependent over $M$. As well, it follows from the `easy half' of Lemma~\ref{nonforking} (not requiring $p$-semiregularity)
that $\clp(a'bM)\cap\clp(cdM)=\clp(M)$.  [In more detail, if there were a `bad' $g$, then taking $g'$ such that $\stp(g'/a'bcdM)=\stp(g/a'bcdM)$ but $\fg {g'} {Ma'bcd} {M^*}$,
we would have $g'\in(\clp(a'bM^*)\cap\clp(cdM^*))\setminus \clp(M^*)$, contradicting $(a',b,c,d)$ a witness to the non-modularity of $p|M^*$ over $M^*$.]
%Thus, $\clp(abM)\cap\clp(cdM)=\clp(M)$ by Lemma~\ref{nonforking}.   That is, $(a',b,c,d)$ is a witness to the non-modularity of $p$ over $M$.
\endproof

The following Corollary codifies what we will use in the subsequent sections.

\begin{Corollary}  \label{whatweuse}
 Suppose $T$ is classifiable, $M\preceq N\models T$,  and $p\in S(M)$ is regular but not locally modular.  If there are $(b,c,d)\in N\setminus M$ realizing
$p^{(3)}$, then there is $a\in N\setminus M$ such that $(a,b,c,d)$ is a witness to the non-modularity of $p$ over $M$.
\end{Corollary}

\bp  Choose an $L(M)$-formula $S$ as in Proposition~\ref{overM}.   Then $\C\models \exists x S(x,b,c,d)$ so $N\models \exists x S(x,b,c,d)$ by elementarity. By Proposition~\ref{overM}, any $a\in N$ satisfying 
$S(a,b,c,d)$ suffices.
\endproof

%\begin{Corollary}  \label{whatweuse}  If $T$ is classifiable, then for any model $M$ and any  regular, non-locally modular type $p\in S(M)$, there
%is an $L(M)$-formula $S(x,y,z,w)$ such that for any $bcd$ realizing $p^{(3)}$, $S(x,b,c,d)$ is consistent and any realization $a$ 
%implies $(a,b,c,d)$ is a witness to the non-modularity of $p$.  
%\end{Corollary}
%
%\bp  All of \ref{4dep}--\ref{overM} were done with the additional assumption that the $L(M)$-formula $\theta(x)$ was 0-definable.  For the general case,
%absorb the parameters $\ebar$ of $\theta$ into the signature.  In the expanded language, obtain $R^*(x,y,z,w,h')$ as in the proof of Corollary~\ref{over M}
%and put  $S(x,y,z,w):=R^*(x,y,z,w,h'\ebar)$.

%and let $e$ be  a finite subset from $M^*$ over which $\tp(abcd/M^*)$ is based and which contains the parameters over which $R^*$ is defined.
%Write $R^*$ as $R^*(x_1,\dots,x_4,e)$.  So $\tp(abce/M)$ implies $\exists x_4 R^*(a,b,c,x_4,e)$.  As $\fg {abc} {M} e$, it follows from finite
%satisfiability that there is some $e'\in M$ with $\exists x_4 R^*(a,b,c,x_4,e')$.  Choose any $d$ witnessing this and check that
%$(a,b,c,d)$ is a witness to the non-modularity of $p$ over $M$.
%\endproof
%
                                                                                                                                                                                                                                                                                                                                                                     
\subsection{Depth-zero like types and minimality}

%The following Lemma is routine, but it is interesting that its proof does not require NDOP (although the verification that a non-trivial
%regular type has depth zero does).  
Recall the definition of a regular type $p$ being of depth zero given in Definition~\ref{depth0def}.  
By working with a-models, it  is easily checked that a regular type being of depth zero is preserved under non-orthogonality.

\begin{Lemma} \label{onedepth0} Suppose $T$ is superstable. 
\begin{enumerate}
\item   For $M$ any model, if $p=\tp(c/M)$ is regular of depth zero, then for any model $M'$ dominated by $c$ over $M$, every  regular $q$ non-orthogonal to $M'$ 
is non-orthogonal to $M$.
\item  Suppose $M\subseteq_{na} N$, $c\in N\setminus M$ is such that $p=\tp(c/M)$ is regular of depth zero, and let $M'\preceq N$ be any model dominated by $c$ over $M$.
Then  $M'$ is minimal over $Mc$ and $M'\subseteq_{na} N$.
\end{enumerate}
\end{Lemma}

\bp  (1)  Fix $M$, $c$, $M'$ and $q$, and suppose $q\in S(A)$.  By superstability choose a finite $B\subseteq M'$ such that $q\not\perp B$.
Choose 
any a-model $M^*\succeq M$ with $\fg {M^*} M {ABc}$.  By Fact~\ref{domfact}(2), $Bc$ is dominated by $c$ over $M^*$, so by Fact~\ref{domfact}(3)
there is an $a$-prime model $M^*[c]$ containing $B$.
As $q\not\perp M^*[c]$ and since $M^*[c]$ is an $a$-model, there is a regular $q'\in S(M^*[c])$ with $q'\not\perp q$.
 As $tp(c/M^*)$ also has depth zero, $q'$ (and hence $q$) is non-orthogonal to $M^*$.  
 As $A$ is independent from $M^*$ over $M$,
we conclude that $q\not\perp M$.

(2)  We first show that $M'$ is minimal over $Mc$.  By way of contradiction, assume there were a proper $M''\prec M'$ containing $Mc$.
By Fact~\ref{regularpair} there is a regular type $q\in S(M'')$ realized in $M'$.   By (1), $q\not\perp M$, so
 as $M\subseteq_{na} M'$, it follows from the 3-model Lemma (Fact~\ref{original3model}) that there is some $a\in M'\setminus M$ with $\fg a M {M''}$.
But this contradicts $M''$ being dominated by $c$ over $M$.  Thus, $M'$ is minimal over $Mc$.  To see that $M'\subseteq_{na} N$,  as $M\subseteq_{na} N$,
apply Fact~\ref{existdominated}, taking $A$ to be $M'$, to obtain $M^*\subseteq_{na} N$ containing $M'$ with $M^*$ dominated by $M'$ over $M$.  As $M'$ was dominated by
$c$ over $M$, the same holds for $M^*$ by Fact~\ref{domfact}(1).    From  the previous sentences, $M^*$ is minimal over $Mc$, hence $M^*=M'$.
\endproof

%Note that in the above, for an arbitrary model $M$, a model $M'$ dominated by $c$ over $M$ need not exist.  However, if $M\subseteq_{na} N$, then we can say more.
The following Lemma extends this  to  independent tuples of depth zero types.

\begin{Lemma} \label{3model}  Suppose $T$ is superstable and $M\na N$.  For every $n$, if $\{c_i:i<n\}\subseteq N$ are $M$-independent with $\tp(c_i/M)$ 
regular of depth zero for each $i$, then a model  $M^*\preceq N$ containing $M\cup\{c_i:i<n\}$ and dominated by $\{c_i:i<n\}$ over $M$ exists.  Moreover,
any such $M^*$ satisfies the following properties:
\begin{enumerate}
\item  $M^*$ is minimal over $M\cup\{c_i:i<n\}$;
\item  $M^*\na N$; 
\item  If $q$ is regular and $q\not\perp M^*$, then $q\not\perp M$.
\end{enumerate}
\end{Lemma}

\bp  For any $n$, the existence of a model $M^*\preceq N$ dominated by $c_0,\dots,c_{n-1}$ over $M$ is immediate by Fact~\ref{existdominated}.
We argue by induction on $n$  that  (1)--(3) hold for any model $M^*\preceq N$ dominated by $c_0,\dots,c_{n-1}$ over $M$.
For $n=1$, this is given by Lemma~\ref{onedepth0}.  
Assume these conclusions  hold for sets of size $n$, and choose $c_0\dots,c_n$ from $N$ satisfying the hypotheses.  Let $M^*$ be dominated by
$\{c_i:i\le n\}$ over $M$.  We first show $M^*$ is minimal over $M\cup\{c_i:i\le n\}$.
Choose any $M'\preceq M^*$ containing $M\cup\{c_i:i\le n\}$.
As $M\subseteq_{na} M'$ we work inside $M'$.
Let $M_0\preceq M'$
be dominated by $\{c_i:i<n\}$ over $M$.
By our inductive hypothesis, both  $M_0\subseteq_{na} M'$ and  any regular type $q$ non-orthogonal to $M_0$ is non-orthogonal to $M$.
By Fact~\ref{existdominated},
 choose $M_1\preceq M'$ to be dominated by $c_n$ over $M_0$.  By Lemma~\ref{onedepth0},  any regular type non-orthogonal to $M_1$ is non-orthogonal to $M_0$.
We argue that $M_1=M^*$, which implies our desired $M'=M^*$.
For this, suppose $M_1\neq M^*$.   By Fact~\ref{regularpair}, there is some $d\in M^*\setminus M_1$ with $q:=\tp(d/M_1)$ regular.
From above, $q$ is also non-orthogonal to $M_0$ and hence to $M$.  
So, by the 3-model Lemma (Fact~\ref{original3model}) applied to the triple $(M,M_1,M^*)$, there would be $e\in M^*$ such that
$\tp(e/M)$ is regular and non-orthogonal to $q$, with $\fg e M {M_1}$.  As $\{c_i:i\le n\}$, this implies $\fg e M {c_0,\dots,c_n}$, contradicting
$M^*$ being dominated by $\{c_i:i\le n\}$ over $M$.
Thus, $M^*$ is minimal over $M\cup\{c_i:i\le n\}$, proving (1).

%  there would be a regular type $q=\tp(d/M_1)$ realized in the difference.  As $\tp(c_n/M_0)$
%is a non-forking extension of $\tp(c_n/M)$, it is regular and of depth zero.  Thus, by definition of depth zero, $q\not\perp M_0$.  Thus, by the inductive
%hypothesis, $q$ is non-orthogonal to $M$.  So, by the 3-model Lemma (Fact~\ref{original3model}) applied to the triple $(M,M_1,M^*)$, there would be $e\in M^*$ such that
%$\tp(e/M)$ is regular and non-orthogonal to $q$, with $\fg e M {M_1}$.  As $\{c_i:i\le n\}$, this implies $\fg e M {c_0,\dots,c_n}$, which contradicts 
%$M^*$ being dominated by $\{c_i:i\le n\}$ over $M$.
%Thus, $M^*$ is minimal over $M_\{c_i:i\le n\}$, proving (1).  

Next, by Fact~\ref{existdominated} there is a model $N'\na N$ that contains and is dominated by $M^*$ over $M$.  By Fact~\ref{domfact}(1),  $N'$ is dominated by
$\{c_i:i\le n\}$ over $M$, so by (1)  applied to $N'$, it is minimal over $M\cup\{c_i:i\le n\}$.
  Thus, $N'=M^*$, so $M^*\na N$, giving (2).  

Finally, choose any regular $q$ non-orthogonal to $M^*=M_1$, using the notation in (1).  Since  $\tp(c_n/M_0)$ has depth zero, $q$ is non-orthogonal to $M_0$ by Lemma~\ref{onedepth0}, hence is non-orthogonal 
to $M$ by our inductive hypothesis.
\endproof

The following definition extends the concept of depth zero to both finite, independent tuples of depth zero types as well as to types dominated by such tuples.

\begin{Definition} \label{depth0like}  {\em  A strong type $p$ is {\em depth-zero like} if every regular type $q$ non-orthogonal to $p$ is of depth zero.
}
\end{Definition}

As examples, if a regular type $p$ has depth zero, then any $p$-semiregular type $q$ is depth zero-like.  
The following Proposition uses classifiability to obviate the  need for  $M\na\C$ in Lemma~\ref{3model}(1).

\begin{Proposition} \label{regmin}  Suppose $p$ is depth zero-like.
\begin{enumerate}
\item  If $T$ is superstable and $M$ is an a-model  on which $p$ is based, then for any realization $b$ of $p|M$ and  every model $N\supseteq Mb$ that is dominated by $b$ over $M$ 
we have:  
\begin{enumerate}
\item   every regular type $q$ non-orthogonal to $N$ is non-orthogonal to $M$: and
\item  $N$ is minimal over $Mb$.
\end{enumerate}
\item  If $T$ is classifiable, then (1) holds for every model $M$ on which $p$ is based.
\end{enumerate}
\end{Proposition}

\bp  (1)  Assume $T$ is superstable.  Fix an a-model $M$ on which $p$ is based and a realization $b$ of $p|M$.  We first show (a) holds when $N=M[b]$, any $a$-prime model over
$Mb$.  
%We first argue that if a regular type $q$ is non-orthogonal to any  a-prime model $M[b]$ over $Mb$, then $q$ is non-orthogonal to $M$.
To see this, choose a maximal $M$-independent set $\{c_1,\dots,c_n\}\subseteq M[b]$
 with $\tp(c_i/M)$ regular.  ($n$ is finite because $\tp(b/M)$ has finite weight.)  As each $c_i$ forks with $b$ over $M$, $\tp(c_i/M)$ is of depth zero.
 Let $M[c_1,\dots,c_n]\preceq M[b]$ be any model dominated by $\{c_1,\dots,c_n)$ over $M$.  We claim that $M[c_1,\dots,c_n]=M[b]$.  For, if not, 
 then by Fact~\ref{regularpair} there would be some $d\in M[b]\setminus M[c_1,\dots,c_n]$ with $r:=\tp(d/M[c_1,\dots,c_n]$ regular.
 By Lemma~\ref{3model}(3), we would have 
 $r\not\perp M$, hence by Fact~\ref{original3model} there would be $e\in M[b]\setminus M$ with $\tp(e/M)$ regular, but $\fg e M {c_1,\dots,c_n}$, contradicting the maximality of
 $\{c_1,\dots,c_n\}$.
Thus,  (a) holds for $M[b]=M[c_1,\dots,c_n]$ by Lemma~\ref{3model}(3).  
For the general case,  take $N\supseteq Mb$ to be any model dominated by $b$ over $M$.  Choose any regular $q$ non-orthogonal to $N$.
By superstability, choose a finite $d\subseteq N$ such that $q\not\perp Mbd$.  As $bd$ is dominated by $b$ over $M$, apply Fact~\ref{domfact}(3) to find an a-prime $M[b]$ over $Mb$
containing $d$.  Then $q\not perp M[b]$, hence $q\not\perp M$ by the special case above.

For (b), choose any model $N\supseteq Mb$ that is dominated by $b$ over $M$.  
Choose any $N'\preceq N$ containing $Mb$ and assume by way of contradiction that $N'\neq N$.  By Fact~\ref{regularpair}, choose $c\in N\setminus N'$ such that $q=\tp(c/N')$ is regular.  By (a) applied to $N'$ and $M$ we have
$q\not\perp M$.
So, by the  3-model Lemma  (Fact~\ref{original3model}) , there is $c^*\in N-M$ that does not fork with $N'$ over $M$.  As $b\in N'$ we conclude that $\fg {c^*} M b$, which contradicts $N$ being dominated by $b$ over $M$.

(2)  Now assume that $T$ is classifiable.
Fix any model $M$ on which $p$ is based and fix a realization $b$ of $p|M$.  
Suppose $N$ is dominated by $b$ over $M$.  To show that $N$ is minimal over $Mb$,
choose any $N_0\preceq N$ containing $Mb$.  To see that $N_0=N$, choose any $c\in N$ and we will conclude that $c\in N_0$. 
To start, choose any a-model $M^*\succeq M$ with $\fg b M {M^*}$.  By Fact~\ref{domfact}(2), $N$ is dominated by $b$ over $M^*$.
By PMOP, choose a constructible  model  $N_2$  over $N\cup M^*$.  By Fact~\ref{lisodom}(2) $N_2$ is dominated by $N$ over $M^*$,
hence $N_2$ is also dominated by $b$ over $M^*$ by Fact~\ref{domfact}(1).   It follows from (1) that
$N_2$ is minimal over $M^*b$.  Also, $c\in N_2$. 

As $N_0\preceq N$, by PMOP again we can find $N_1\preceq N_2$ that is constructible over $N_0\cup M^*$.  The minimality of $N_2$ over $M^*b$ implies that
$N_2=N_1$, hence $c\in N_1$.  As $N_1$ is atomic over $N_0M^*$, we have that $\tp(c/N_0M^*)$ is isolated.  However, as $cN_0$ is dominated by $b$ over $M$,
the fact that $\fg b M {M^*}$ implies that $\fg {cN_0} M {M^*}$.    As $M\preceq N_0$, we have $\fg c {N_0} {M^*}$.  
The  Open Mapping Theorem implies that $\tp(c/N_0)$ is isolated, hence $c\in N_0$.
\endproof

\subsection{The main theorem}

Lemma~\ref{3model} suggests the following notation.
 If $M\na\C$ and $\{c_i:i\le n\}$ are $M$-independent and realize regular, depth zero types over $M$, then
the notation $M(c_0,\dots,c_n)$ refers to any model dominated by $c_0\dots,c_n$ over $M$.  The reader is cautioned that even when $\tp(b/M)=\tp(c/M)$
are the same regular type, the notations $M(b)$ and $M(c)$ represent possibly different isomorphism types of models dominated over $M$
by $b$ or $c$, respectively.  

Now suppose that $T$ is classifiable, so in particular, a constructible model $N$ exists over any independent triple of models.

\begin{Lemma}  \label{combine4}  Suppose $T$ is classifiable, $\{c_i:i<n\}$ independent over a model $M$ and for each $i$, $\tp(c_i/M)$ has depth zero and
$M(c_i)$ is dominated by $c_i$ over $M$.
Then the constructible model $N$ over $\bigcup\{M(c_i):i<n\}$ is dominated by $\{c_i:i<n\}$ over $M$.  Thus, $N$ is minimal over $M\cup\{c_i:i<n\}$ and its isomorphism type over $M$ is uniquely determined by the isomorphism types $\{M(c_i):i<n\}$ over $M$.  
\end{Lemma}

\bp   As $\{c_i:i<n\}$ are $M$-independent and as each $M(c_i)$ is dominated by $c_i$ over $M$, it follows by forking calculus that $\{M(c_i):i<n\}$ are $M$-independent.
Thus,  by PMOP, a constructible model $N$ over $A:=\bigcup\{M_i:i<n\}$ exists.  By Fact~\ref{lisodom}(2) $N$ is dominated by $A$ over $M$.  Thus,
$N$ is dominated by $\{c_i:i<n\}$ over $M$ by Fact~\ref{domfact}(1).  By Lemma~\ref{3model}(1), $N$ is minimal over $M\cup\{c_i:i<n\}$.   As $N$ is constructible over $A$, its isomorphism type over $M$
depends only on the isomorphism type of $A$ over $M$.  By independence and $M$ being a model, the isomorphism type over $A$ over $M$ is determined by the set of isomorphism types of its component pieces.
\endproof

In light of Lemma~\ref{combine4}, we write such an $N$ as $M(c_0,\dots, c_{n-1})$.     
In particular, if $\tp(c/M)$ has depth zero and $M(c)$ is chosen, we let $M(c)^{(4)}$ denote the prime model over four $M$-independent copies of $M(c)$.
It follows from Lemma~\ref{combine4} that  the isomorphism type of
$M(c)^{(4)}$ over $M$ is uniquely determined by the type of $M(c)$ over $M$.

%Proposition~\ref{main} is the most technical result of this paper. 
 For the next few results, recall that by Fact~\ref{NDOP}, any non-locally modular type has depth zero.
 Proposition~\ref{main} is the most technical result of this paper. 

\begin{Proposition}  \label{main}
Suppose $T$ is classifiable, $M\na \C$, $p\in S(M)$ is regular but not locally modular.
  Then for any realizations $b,c$ of $p$ and any choice of  models $M(b)$ and $M(c)$ as above, 
$M(b)$  elementarily embeds into $M(c)^{(4)}$ over $M$.
\end{Proposition}

\bp  
Without loss of generality, we may assume $b$ and $c$ are independent.
Fix choices for $M(b)$ and $M(c)$.  Choose $d$ realizing $p|Mbc$ and choose $M(d)$ to be isomorphic to $M(c)$ over $M$.
Let $M(bcd)$ denote the prime model over $M(b)M(c)M(d)$.  Since $p$ is not locally modular 
and as $(b,c,d)$ realizes $p^{(3)}$, Corollary~\ref{whatweuse} gives an $a\in M(bcd)$ such that $(a,b,c,d)$ is a  witness to non-modularity
over $M$. Since $M\na M(bcd)$, we can use Fact~\ref{existdominated} to find some $M(a)\subseteq_{na} M(bcd)$.  
As $M(a),M(b),M(c)\preceq M(bcd)$, by PMOP there is a prime model $M^*\preceq M(bcd)$ over $M(a)M(b)M(c)$.  
Clearly, $M^*$ is dominated by $abc$ over $M$ and $(a,b,c)$ realizes $p^{(3)}$, so we can write it as $M(abc)$.

\medskip\noindent{\bf Claim 1.}  $M^*=M(bcd)$.

\bp  If not, then use Fact~\ref{regularpair} to choose $g\in M(bcd)\setminus M^*$ to realize a regular type $q$.  
By Lemma~\ref{3model} we have $q\not\perp M$, so by the 3-model Lemma (Fact \ref{original3model}) there would be $h\in M(bcd)$ such that $\tp(h/M)$ is regular and
non-orthogonal to $q$, with $\fg h M {M^*}$.  This $h$ is dominated by $bcd$ over $M$, so $\tp(h/M)$ must be non-orthogonal to $p$.
But then ${h,a,b,c}$ are independent realizations of regular types non-orthogonal to $p$, contradicting $w_p(M(bcd)/M)=3$.
\endproof

So $M^*$ is equal to both $M(abc)$ and $M(bcd)$.  Next, let $M(ab)\preceq M(abc)$ be the unique model that is prime over $M(a)\cup M(b)$.
For the moment we work over $M(ab)$.  Choose $e$ to realize $p|Mabcd$  and fix an isomorphism 
$$\Phi:M(abc)\rightarrow M(abe)$$ 
fixing $M(ab)$ pointwise with $\Phi(c)=e$.  Put $M(e):=\Phi(M(c))$.
As both $d$ and $M(d)$ are contained in $M(abc)$, we let $f:=\Phi(d)$ and $M(f):=\Phi(M(d))$.
Let $N^*$ be prime over $M(abc)\cup M(abe)$ over $M(ab)$.  Note that $N^*$ is dominated by $ce$ over
$M(ab)$.  

We now use the fact that $(a,b,c,d)$ is a witness to non-modularity.
By construction, $ef$ and $cd$ have the same type over $M(ab)$, so it follows from Lemma~\ref{firstL} that $\{c,d,e,f\}$
are independent over $M$.  It follows that the four models $M(c),M(d),M(e),M(f)$ are $M$-independent.  Since $M(d)$ was chosen to be isomorphic to $M(c)$ over $M$,
 $M(e)=\Phi(M(c))$,  and $M(f)=\Phi(M(d))$, the four models are pairwise isomorphic over $M$.

As each of $M(c),M(d),M(e),M(f)$ are contained in $N^*$, let
$N\preceq N^*$ be constructible over $M(c)\cup M(d)\cup M(e)\cup M(f)$.  Note that $N$ is isomorphic to $M(c)^{(4)}$ over $M$
and $M(b)\preceq M(ab)\preceq N^*$.
Thus, the Proposition is proved once we establish the following claim.

\medskip
\noindent
{\bf Claim 2.}  $N=N^*$.

\bp  If not, then by Fact~\ref{regularpair} choose $g\in N^*\setminus N$ with $q:=\tp(g/N)$ regular.  
By Lemma~\ref{3model}, $q$ is non-orthogonal to $M$.  Thus, by the 3-model Lemma (Fact~\ref{original3model}),
there is $h\in N^*\setminus M$ such that $\tp(h/M)$ is regular and non-orthogonal to $q$, with $\fg h M N$.
We split into cases depending on the non-orthogonality class of $q$.

First, assume that $q$ is non-orthogonal to $p$.
On one hand, $\{h,c,d,e,f\}\subseteq N^*$ consists of 5 independent realizations of regular types non-orthogonal to $p$.
On the other hand, $w_p(M(ab)/M)=2$ and $w_p(N^*/M(ab))=2$, so $w_p(N^*/M)=4$, which is a contradiction.

Finally, assume that $q$ is orthogonal to $p$.  Then clearly $\tp(h/M(ab))$ does not fork over $M$.  But $N^*$ is dominated by
$ce$ over $M(ab)$, and by the orthogonality, $\fg h {M(ab)} {ce}$.  Thus, by transitivity, $\fg h M N^*$, which is absurd since
$h\in N^*\setminus M$.
\endproof

\begin{Corollary}  \label{manytp}  $T$ classifiable.  Suppose $M\na\C$ is countable and $p:=\tp(c/M)$ is regular but not locally modular.
Then for each $n$, $Q_n:=\{q\in S_n(Mc): \hbox{$\dbar c$ is dominated by $c$ over $M$ when $\dbar \models q$}\}$ is countable.
\end{Corollary}

\bp  Fix a model $M(c)$  dominated by $c$ over $M$ and fix an $n\ge 1$.  
For any $q\in Q_n$, choose a realization $\dbar_q$ of $q$ and let $N_q$ be any countable, $\ell$-constructible model over $Mc\dbar_q$.
By Fact~\ref{lisodom}(2) $N_q$  is dominated by $\dbar_q c$ (and hence by $c$) over $M$.  By Proposition~\ref{main}, choose an elementary embedding
$f_q:N_q\rightarrow M(c)^{(4)}$ fixing $M$ pointwise.  If $Q_n$ were uncountable, there would be distinct $q\neq q'$ with $f_q(\dbar_q c)=f_{q'}(\dbar_{q'} c)$.
Thus, $\tp(\dbar_qcM)=\tp(\dbar_{q'}cM)$, i.e., $\tp(\dbar_q/Mc)=\tp(\dbar_{q'}/Mc)$,
contradicting $q\neq q'$.
\endproof

The following technical Lemma, which is true in any superstable theory,
 is implicit in the proof of Lemma~5.3 of \cite{ShelahBuechler}, but is included here at the request of the referee. 

\begin{Lemma}  \label{tech53}  $T$ superstable. Suppose $M\subseteq_{na}\C$ and $A\supseteq M$.  For every non-algebraic $\psi(x,a_0)$ (with $a_0\subseteq A$),
suppose $\phi(x,a)$ is of least $R^\infty$-rank satisfying the following conditions:
\begin{itemize}
\item  $a\subseteq A$,  $\phi(x,a)\vdash \psi(x,a_0)$, and $\phi(x,a)$ consistent, but non-algebraic.
\end{itemize}
Then $cA$ is dominated by $A$ over $M$ for every $c$ realizing $\phi(x,a)$.
\end{Lemma}

\bp  The driving force is the equivalent of $M\subseteq_{na} \C$ proffered by Lemma~2.9(c) of 
\cite{HrShOT}.  Choose any such $\phi(x,a)$ and by way of contradiction, assume  that  $cA$ is not dominated by $A$ over $M$ for some $c$ with $\phi(c,a)$.
Choose an element $b\in\C$ such that $\fg b M A$, but $\nfg b M {Ac}$.  Possibly enlarging $a$ within $A$, choose an $L(M)$-formula $\theta(x,y,z)\in\tp(cab/M)$
such that $\theta(x,a,b)$ forks over $Ma$.  Fix a finite $F\subseteq M$ containing the hidden parameters from $\theta$ with $\fg {cab} F M$.
Then $\fg a F {Mb}$ and $\phi(x,a)\wedge\theta(x,a,b)$ is consistent and forks over $Fa$.  Thus, by Lemma~2.9(c) of \cite{HrShOT}, there is some $b'\in M$ and some
$L(M)$-formula $\theta'(x,a,b')$ such that $\phi(x,a)\wedge\theta'(x,a,b')$ is consistent and forks over $Fa$.  
Thus, $R^\infty(\phi(x,a)\wedge \theta(x,a,b'))<R^\infty(\phi(x,a))$.  As $b'\subseteq A$, this  contradicts the minimality of $R^\infty(\phi(x,a))$.
\endproof

\begin{Theorem}  \label{atomic} $T$ classifiable.
Suppose $M\na\C$ is countable and $p:=\tp(c/M)$ is regular but not locally modular.
%$\tp(c/M)=p$, and $p$ has a witness to non-local modularity  over   $M$.
Then there is a constructible, minimal  model over $Mc$.
\end{Theorem}

\bp  Assume by way of contradiction that there is no constructible model over $Mc$.
As $Mc$ is countable, there is a finite $\bbar$
such that $\tp(\bbar/Mc)$ is isolated, and a consistent, non-algebraic formula $\theta(x)$ over $M\bbar c$ that has no complete extension $\psi(x)\vdash\theta(x)$ over $Mc\bbar$.
Among all consistent, non-algebraic formulas $\phi(x,\bbar c)\vdash \theta(x)$, choose one of minimal $R^\infty$-rank.  By Lemma~\ref{tech53} $d\bbar c$ is dominated by
$\bbar c$ over $M$ for every realization of $\phi(x,\bbar c)$.  However, $\tp(\bbar/Mc)$ is isolated, hence $\bbar c$ is dominated by $c$ over $M$ by Fact~\ref{lisodom}(2).
Thus,  by Fact~\ref{domfact}(1), 
$d\bbar c$ is dominated by $c$ over $M$ for every $d$ realizing $\phi(x,\bbar c)$.  However, $\phi(x,\bbar c)$ has no complete extension, so there are a perfect
set of complete types in $S_k(M\bbar c)$ extending $\phi(x,\bbar c)$, where $k=\lg(x)$.  Thus, there are uncountably many distinct types in $Q_n$ (where $n=k+\lg(\bbar)$)
contradicting Corollary~\ref{manytp}.

Thus, a constructible model $N$ over $Mc$ exists.  As $N$ is dominated by $c$ over $M$ by Fact~\ref{lisodom}(2), its minimality over $M$ follows from Lemma~\ref{onedepth0}(2).
\endproof

%Fix such a $\bbar$ and $\theta$.  Choose a  consistent formula $\psi(x)\vdash\theta(x)$, also over $Mc\bbar$, of smallest $R^{\infty}$-rank.  
%By the final sentence of Fact~\ref{existdominated}, since $M\subseteq_{na}\C$,  every realization $d$ of $\psi(x)$ has $d\bbar c$ dominated by $\bbar c$ over $M$.
%As $\tp(\bbar/Mc)$ is isolated, Fact~\ref{lisodom} implies $\bbar c$ is dominated by $c$ over $M$, we conclude that $ d\bbar c$ is dominated by $c$ over $M$.
%
%
%
%However, as there is no complete $\phi\vdash\psi$, there are
% $2^{\aleph_0}$ types in $S_1(Mc\bbar)$ containing $\psi$.    It follows that there are $2^{\aleph_0}$ types over $Mc$ that are dominated by $c$ over $M$, contradicting
% Corollary~\ref{manytp}.    
%
%Thus, a constructible model $N$ over $Mc$ exists.  As $N$ is dominated by $c$ over $M$, its minimality over $M$ follows from Lemma~\ref{onedepth0}(2).
%\endproof

The following strengthening,  obviating the requirements that  $M$ be countable and $M\subseteq_{na} \C$,  is the main result of this section.

\begin{Theorem}  \label{atomic2} Let $M$ be any model of a classifiable theory $T$.
%Suppose that $T$ is classifiable and that
%$M$ is any model. 
If $p\in S(M)$ is regular but not locally modular and $b$ is any realization of $p$ 
then every model $N$  containing $Mb$ that is dominated by $b$ over $M$ is both constructible and minimal over $Mb$.
\end{Theorem}

\bp    As $T$ is classifiable and $p$ is non-locally modular, then $p$ has depth zero by Fact~\ref{NDOP}.
We first prove the Theorem when $M$ is countable.  Choose any $N\succeq M$ that is dominated by $b$ over $M$, e.g., an $\ell$-constructible model over $Mb$.
By  Proposition~\ref{regmin}(2), $N$ is minimal over $Mb$, hence $N$ must be countable as well.

To show that $N$ is constructible over $Mb$,  choose any countable $M^*\succeq M$ with $\fg N M {M^*}$ such that $M^*\na\C$.  
 By PMOP, choose a countable $N^*$ that is prime over $N\cup M^*$.  
% By Theorem~\ref{atomic}, there is $S$, constructible and minimal over $M^*b$.
% As $S$ is prime over $Mb$, we may assume $S  \preceq N^*$.
%We claim that $N^* = S$.  
As $\fg N M {M^*}$ and $N$ is dominated by $b$ over $M$, $N$ is also dominated by $b$ over $M^*$ by Fact~\ref{domfact}(2).
As $N^*$ is atomic over $M^* \cup N$,  it is dominated by $N$ over
$M^*$ by Fact~\ref{lisodom}(2), hence $N^*$ is dominated by $b$ over $M^*$ by Fact~\ref{domfact}(2).
 By Lemma~\ref{onedepth0} $N^*$ is minimal over $M^* b$.   
  By Theorem~\ref{atomic},  as $M^*\na\C$ there is a countable $S$, constructible and minimal over $M^*b$.
 As $S$ is prime over $M^*b$, we may assume $S  \preceq N^*$.
As $N^*$ is minimal over $M^*b$,  $N^* = S$, so $N^*$ is constructible over$M^*b$.
To finish, since $N \subset N^*$, $N$  (as a set)  is atomic over $M^* b$. By the Open mapping theorem, it follows that $N$ is already atomic over $Mb$. As $N$ is a countable, it is constructible over $Mb$.

Now suppose $M$ is arbitrary.  Choose any $N\supseteq Mb$ that is dominated by $b$ over $M$.
By Proposition~\ref{regmin}(2), $N$ is minimal over $Mb$.  
To show that $N$ is constructible over $Mb$, 
consider the pair $(N,M)$ in a language with a predicate $U$ for $M$.  Take  a countable elementary substructure  in the pair language, 
$(N',M')\preceq (N,M)$ with $b\in N'$. It follows that, in the original language, $N'$ and $M$ are independent over $M'$.

As $N' \subset N$, $b$ dominates $N'$ over $M$. By the independence of $N'$ and $M$ over $M'$, $b$ also dominates $N'$ over $M'$ by Fact~\ref{domfact}(2).  
As $N'$ is countable, it follows from the argument above that $N'$ is constructible over
$M'b$ by a construction sequence $\cbar=\<c_\alpha:\alpha<\beta\>$.  By Fact~\ref{freeconst},  $\cbar$ is also a construction sequence over $Mb$.
Now by PMOP,  take $N^*\preceq N$, a constructible model over $N'$ and $M$, inside $N$.  As $Mb\cbar=N'$ as a set, it follows that $N^*$ is constructible over $Mb$.
As $N$ is minimal over $Mb$, $N^* = N$ so $N$ is constructible over $Mb$.   
%To see that this implies that $N$ is atomic over $Mb$ is a routine Tarski-Vaught argument. [Choose a formula $\phi(x,b,m')$
%isolating $\tp(e/M'b)$.  Then, for every $L$-formula $\delta(x,y,\zbar)$
%$$(N',M')\models (\forall \zbar\in M')\forall x[\phi(x,b,m')\  \hbox{decides}\ 
%\delta(x,b,\zbar)]$$
%which transfers up to $(N,M)$ by elementarity.]
\endproof

\begin{Corollary}  \label{weight1}  $T$ classifiable.  Suppose that $M\preceq N$ and $N/M$ has weight one, non-orthogonal to a non-locally modular type $p$.
Then $N$ is both constructible and minimal over $Mb$ for any element $b\in N\setminus M$.
\end{Corollary}

\bp   Choose any $b\in N\setminus M$.  
%Recall that $N/M$ having weight one means that  
%$wt(c/M)=1$ for all finite tuples from $N\setminus M$.  
%In particular, $\nfg c M b$ for any tuple $c\in N\setminus M$.
    Since $M\prec N$, by Fact~\ref{regularpair}, choose $a\in N\setminus M$ such that $p':=\tp(a/M)$ is regular.  As $N/M$ is weight one, non-orthogonal
to $p$, $p'\not\perp p$, so $p'$ is also non-locally modular (see e.g., 7.2.4 of \cite{Pillay}).   By Theorem~\ref{atomic2}, $N$ is constructible and minimal over $Ma$.
In particular, $\tp(b/Ma)$ is isolated
by some $L(M)$-formula $\alpha(a,y)$.  Dually, we have the following.

\medskip
\noindent
{\bf Claim.}  $\tp(a/Mb)$ is isolated.
\noindent
\medskip

\bp  As $p'\in S(M)$ is non-locally modular, choose $\theta(x)\in p'$ as in Theorem~\ref{olddefinability}.  Also, since $N/M$ has weight one, $\nfg a M b$.  Choose an L(M)-formula
$\beta(x,b)$ witnessing this forking.  Put 
$$\psi(x,b):=\theta(x)\wedge\alpha(x,b)\wedge\beta(x,b)$$ 
Clearly, $\psi(a,b)$ holds and we show it isolates $\tp(a/Mb)$.
For this, choose any $L(M)$-formula $\delta(x,y)$ and assume $\delta(a,b)$ holds.  As $N\preceq\C$, it suffices to show that any $a'\in N$ realizing $\psi(x,b)$ satisfies $\delta(x,b)$.
So choose any $a'\in \psi(N,b)$.  Since $\theta(a')$ holds, $\tp(a'/M)$ is $p$-simple.  Since $\beta(a',b)$ holds, $\tp(a'/Mb)$ forks over $M$, so $a'\not\in M$.  Thus, as $N/M$ has weight one, $\nfg a M {a'}$.    Since $\tp(a/M)$ is regular, this implies $w_p(a'/M)>0$, so by Theorem~\ref{olddefinability}(3), $\tp(a'/M)=p'=\tp(a/M)$.
Finally, as $\alpha(a,y)$ isolates $\tp(b/Ma)$, we have $\forall y (\alpha(a,y)\rightarrow \delta(a,y))\in \tp(a/M)$.  Thus, 
$\forall y (\alpha(a',y)\rightarrow \delta(a',y))$ holds as well.  Since $\alpha(a',b)$ was assumed to hold,  we conclude that $\delta(a',b)$ holds..
\endproof

\medskip

Given the Claim,  $N$ is constructible over $Mb$, since from above it  is constructible over $Mab$ and $\tp(a/Mb)$ is isolated.
  To see that $N$ is minimal over $Mb$, note that in fact, $N$ is minimal
over $Ma'$ for any realization of $p'$ in $N$.   So, choose any model $N'\preceq N$ containing $Mb$.  By the Claim there is some $a'\in N'$ realizing $p'$, hence $N'=N$.
\endproof

%As $b/M$ is not orthogonal to $p$, there is $a\in\dcl(Mb)\setminus M$ that is $p$-simple of positive $p$-weight (Fact~\ref{psimpleindcl}).
%As $a\in N$ and as $N/M$ has weight one, it follows that $\tp(a/M)$ is regular, non-orthogonal to $p$.  Thus, by Theorem~\ref{atomic2}, $N$ is constructible and minimal
%over $Ma$.  As $a\in\dcl(Mb)$, it follows that $N$ is constructible and minimal over $Mb$ as well.  
%\endproof

\section{When domination implies isolation}

We begin this section with a recasting of Theorem~\ref{atomic2}.

\begin{Corollary} \label{dom}   Suppose $T$ is classifiable,  $A$ is any set, $p\in S(A)$ is a regular, stationary, non-locally modular type, and $b$ is any realization of $p$.
Then for any $e$, if $be$ is dominated by $b$ over $A$, $\tp(e/Ab)$ is isolated.
\end{Corollary}

\bp  Let $M\na\C$ be free from $b$ over $A$ and choose an $\ell$-constructible model $N\supseteq Mbe$ over $Mbe$.
As $be$ is dominated by $b$ over $M$, hence $N$ is dominated by $b$ over $M$.
By Theorem~\ref{atomic2}, $N$ is constructible, hence atomic over $Mb$, so $\tp(e/Mb)$ is isolated.
Also, $\fg e {Ab} {Mb}$, hence $\tp(e/Ab)$ is isolated by the Open Mapping Theorem.
\endproof

This result suggests the following definition.

\begin{Definition}  \label{DIdef}  {\em
A strong type $p$ satisfies {\em DI} (read `domination implies isolation') if, for every set $A$ on which $p$ is based and stationary and for every realization $b$ of $p|A$,
for every $c\in\C$, if $bc$ is dominated by $b$ over $A$, then $\tp(c/Ab)$ is isolated.
}
\end{Definition}  

In the remainder of this section, we explore this notion in classifiable theories.
%We begin with two results that only require stability.
%
%
%
%
%
%\begin{Lemma}  \label{modisol}  Suppose $T$ is stable,
% $M$ is a model,  $A\supseteq M$  is any set, and $\phi(x,a)$ isolates a type $p\in S(A)$  Then:
%\begin{enumerate}
%\item  For every $\fg B M A$, $p$ has a unique extension $q\in S(AB)$, and $q$ is also isolated by $\phi(x,a)$; and
%\item  For any $c$ realizing $\phi(x,a)$, $cA$ is dominated by $A$ over $M$.
%\end{enumerate}
%\end{Lemma}  
%
%\bp  
%
%For (1), the unique type  is $q:=\{\psi(x,b,a'):\phi(x,a)\vdash d_r y \psi(x,y,a')$, where $a'\in A$, $b\in B$, $\psi(x,y,z)$ over $M$, and $r=\tp(b/M)\}$.
%
%For (2), choose any $B$ satisfying $\fg B M A$.  As $p$ has a non-forking extension to $S(AB)$, it must be the $q$ from (1). 
%\endproof
%
%
Among depth zero-like types, the notion of DI has many equivalents.

\begin{Proposition}  \label{DIprop}  
Suppose $T$ is classifiable. The following are equivalent for a depth zero-like strong type $p$:
\begin{enumerate}
\item  $p$ is DI;
\item  For every countable $M$ on which $p$ is based and for every $b$ realizing $p|M$, and for every $n$, the isolated types in $S_n(Mb)$ are dense;
\item  For every countable $M$ on which $p$ is based, for every $b$ realizing $p|M$, there is a constructible model $N$ over $Mb$.  Moreover,
every model $N$ that is dominated by $b$ over $M$ is constructible over $Mb$;
\item  Same as (3), but for every model $M$ on which $p$ is based;
\item  There is some a-model $M$ on which $p$ is based and some $b$ realizing $p|M$ for which there is a constructible model $N$ over $Mb$.
\end{enumerate}
\end{Proposition}

\bp  We prove $(1)\Rightarrow(2)\Rightarrow(3)\Rightarrow(4)\Rightarrow(1)$.  Then $(4)\Rightarrow(5)$ is trivial and we will show $(5)\Rightarrow(2)$.

$(1)\Rightarrow(2)$.  Assume $(1)$ and choose a countable $M$, $b$, and $n$ as in (2).  Let $\phi(x,b,m)$ be any consistent formula with $\lg(x)=n$.
Choose any model $N$ that is dominated by $b$ over $M$ (e.g., an $\ell$-constructible one).   
Choose any $c\in N$ realizing $\phi(x,b,m)$.  As $cb$ is dominated by $b$ over $M$,
it follows from (1) that $\tp(c/Mb)$ is isolated.

$(2)\Rightarrow(3)$.  Fix $M$ and $b$ as in (3).  As $M$ is countable and the isolated types are dense, it follows from Vaught that a constructible model $N$ over $Mb$
exists.  For the final sentence, let $N^*$ be any model dominated by $b$ over $M$.  As $N$ is prime over $Mb$, we may assume $N\preceq N^*$.  But, as $N^*$ is minimal over
$Mb$ by Proposition~\ref{regmin}(2), $N=N^*$, so $N^*$ is constructible over $Mb$.

$(3)\Rightarrow(4)$.  Fix $M$ and $b$ as in (4).  Choose a countable model $M_0\preceq M$ such that $\fg b {M_0} M$.  By (3), let $N_0$ be constructible over $M_0b$.
Fix a construction sequence $\cbar=\<c_\alpha:\alpha<\beta\>$ for $N_0$ over $M_0b$.  By Fact~\ref{freeconst} $\cbar$ is a construction sequence over $Mb$.
Now, by PMOP, let $N$ be constructible over $N_0M$.  By concatenation, it follows that $N$ is constructible over $Mb$.  For the final sentence, let $N^*$ be any model dominated by
$b$ over $M$.  As $N$ is prime over $Mb$, we may assume $N\preceq N^*$.  But then $N=N^*$, again by Proposition~\ref{regmin}(2).

$(4)\Rightarrow(1)$.  Choose any set $A$ on which $p$ is based and stationary and let $b$ be any realization of $p|A$.  Choose any element $c\in \C$
such that $bc$ is dominated by $b$ over $A$.  Choose any model $M\supseteq A$ satisfying
$\fg M A {bc}$.  It follows that $bc$ is dominated by $b$ over $M$.  Choose any model $N$ dominated by $bc$ over $M$.  
By transitivity we have that $N$ is dominated by $b$ over $M$.  Thus, by (4), $N$ is constructible, and hence atomic over $Mb$.  In particular, $\tp(c/Mb)$ is isolated.
Thus, $\tp(c/Ab)$ is also isolated by the Open Mapping Theorem. 

$(4)\Rightarrow(5)$ is immediate.  

$(5)\Rightarrow(2)$.  Let $M^*$ and $b^*$ be any a-model and witness exemplifying (5).  
Given any countable $M$ on which $p$ is based and $b$ realizing $p|M$, the saturation of $M^*$ implies there is an $M'\preceq M^*$ such that
$\tp(Mb)=\tp(M'b^*)$.  Thus, without loss, assume $M\preceq M^*$ and $b=b^*$.  As $b$ realizes $p|M^*$, we have that $\fg b {M} {M^*}$.
Let $\phi(x,b,m)$ be any consistent formula with $m\in M$.  
By (5), choose a constructible (and hence atomic) $N$ over $M^*b$.  Choose any $c\in N$ realizing $\phi(x,b,m)$.  
Then $\tp(c/M^*b)$ is isolated, hence $\tp(c/Mb)$ is as well by the Open Mapping Theorem.
\endproof

\begin{Proposition}  \label{persevere}
  Let $T$ be classifiable.
  \begin{enumerate}
\item  Suppose $\stp(b/A)$ is depth zero-like and DI and $bc$ is dominated by $b$ over $A$.  Then $\stp(bc/A)$ is both depth zero-like and DI.
\item  Suppose $\tp(bc/M)$ is depth zero-like and DI, and $\tp(c/Mb)$ is isolated.  Then $\tp(b/M)$ is depth zero-like and DI as well.
\item  If $p$ and $q$ are both depth zero-like and DI, then so is $p\otimes q$.
\end{enumerate}
\end{Proposition}

%\bp  (1)  That $\stp(bc/A)$ is depth zero-like is clear.  As for DI, choose an a-model $M$ and $b$ as in Proposition~\ref{DIprop}(5).  Without loss, we may assume
%$A\subseteq M$, so $\fg b A M$.  It follows that $bc$ is dominated by $b$ over $M$, hence $\tp(c/Mb)$ is a-isolated.  Thus, there is an a-prime model $M[b]$ over
%$Mb$ with $c\in M[b]$.  [Think of $c$ as being the first step of an a-construction sequence over $Mb$.]  Then $M[b]$ is also a-prime over $Mbc$.  However, $M[b]$ is
%constructible over $Mb$, so $\tp(c/Mb)$ is actually isolated.  It follows that $M[b]$ is constructible over $Mbc$, so $\stp(bc/A)$ is DI by Proposition~\ref{DIprop}(5).

\bp  (1)  That $\stp(bc/A)$ is depth zero-like is clear.  As for DI, choose an a-model $M\supseteq A$ with $\fg b A M$.  By Fact~\ref{domfact} $bc$ is dominated by $b$ over $M$, hence there is an $a$-prime model $M[b]$ over $Mb$ with $c\in M[b]$.  As $\tp(b/M)$ is DI, $M[b]$ is constructible over $Mb$ by Proposition~\ref{DIprop}(4).  Thus, $\tp(c/Mb)$ is isolated,
so $M[b]$ is also constructible over $Mbc$.  Hence $\stp(bc/A)$ is DI by Proposition~\ref{DIprop}(5).
%and $b$ as in Proposition~\ref{DIprop}(5).  Without loss, we may assume
%$A\subseteq M$, so $\fg b A M$.  It follows that $bc$ is dominated by $b$ over $M$, hence $\tp(c/Mb)$ is a-isolated.  Thus, there is an a-prime model $M[b]$ over
%$Mb$ with $c\in M[b]$.  [Think of $c$ as being the first step of an a-construction sequence over $Mb$.]  Then $M[b]$ is also a-prime over $Mbc$.  However, $M[b]$ is
%constructible over $Mb$, so $\tp(c/Mb)$ is actually isolated.  It follows that $M[b]$ is constructible over $Mbc$, so $\stp(bc/A)$ is DI by Proposition~\ref{DIprop}(5).

(2)  That $\tp(b/M)$ is depth zero-like is immediate.  For DI,  choose an a-model $M^*$ as in Proposition~\ref{DIprop}(5) 
witnessing that $\tp(bc/M)$ is DI.  Without loss, we may assume $M\preceq M^*$,
so $\fg {bc} M {M^*}$.  It follows from Lemma~\ref{lisodom}(3) that $\tp(c/M^*b)$ is isolated as well.  As $\tp(bc/M^*)$ is DI, let $N^*$ be constructible over $M^*bc$.
As $\tp(c/M^*b)$ is isolated, $N^*$ is also constructible over $M^*b$.  Thus, $\tp(b/M^*)$ is DI by Proposition~\ref{DIprop}(5).

(3)  As both $p$ and $q$ are depth zero-like, it is immediate that $p\otimes q$ is as well.
As for DI, let $M$ be any model on which $p\otimes q$ is based and let $(c_1,c_2)$ realize $p\otimes q$.  
As both $p$ and $q$ are DI, there is a constructible model $N_1$ over $Mc_1$ and a constructible model $N_2$ over $Mc_2$.  
As $\fg {c_1} M {c_2}$, domination implies $\fg {N_1} M {N_2}$.  It follows from Fact~\ref{freeconst} that the set $N_2$ is constructible over $Mc_1c_2$.
By Fact~\ref{freeconst} again, the set $N_1$ is constructible over $N_2c_1$.  Concatenating these construction sequences, 
we get that the set $N_1\cup N_2$ is constructible over $Mc_1c_2$.
%by iterating Lemma~\ref{modisol} it follows that $N_1$ is constructible over $Mc_1c_2$ and $N_2$ is constructible over $N_1c_2$.
Finally,  by PMOP there is a model $N^*$ that is constructible over $N_1N_2$.  It follows that $N^*$ is constructible over $Mc_1c_2$, so 
$p\otimes q$ is DI by Proposition~\ref{DIprop}(4).
\endproof

We can combine several of our results in  the following Corollary which  both generalizes Corollary~\ref{weight1} and extends Theorem~\ref{strong}.

\begin{Corollary} \label{end} $T$ classifiable.  Suppose $p$ is a regular, non-trivial, depth zero, DI type (for example a non-locally modular regular type) and let $\tp(a/M)$ be $p$-semiregular of weight $k$.  Then there is a constructible, minimal model $N$ over $Ma$.
\end{Corollary}

\bp  Choose an $a$-model $M^*$ independent from $a$ over $M$ and choose an $M^*$-independent tuple $\cbar=\<c_i:i<k\>$ of realizations of $p|M^*$ such that $a$ and $\cbar$
are domination equivalent over $M^*$. 

As $\tp(a/M^*)$ is $p$-semiregular, it is depth-zero like. By Proposition~\ref{persevere}(3) and (1),  both $\tp(\cbar/M^*)$ 
and $tp(a\cbar/M^*)$ are depth zero-like and $DI$.  To show that $\tp(a/M^*)$ is DI, by Proposition~\ref{persevere}(2) it suffices to show $\tp(\cbar/M^*a)$ is isolated.

To see the isolation, first note that since $\tp(\cbar/M^*)$ is DI, there is an $L(M^*)$-formula $\alpha(x,\cbar)$ isolating $\tp(a/M^*\cbar)$.
Next, choose  $L(M^*)$-formulas  $\phi(x)\in\tp(a/M^*)$ and $\theta(y)\in p|M^*$ as in Fact~\ref{pweightprov} and Theorem~\ref{olddefinability}, respectively.
As $w_p(a/M^*\cbar)=0$,  the definability of $p$-weight  inside $\phi(x)$
implies there is a formula $\beta(a,\ybar)\in\tp(\cbar/M^*a)$ such that $w_p(a/M^*\dbar)=0$ for all $\dbar$ realizing
$\beta(a,\ybar)$.
Put 
$$\psi(a,\ybar):=\bigwedge_{i<k} \theta(y_i)\wedge \alpha(a,\ybar)\wedge\beta(a,\ybar)$$
Visibly, $\psi(a,\ybar)\in\tp(\cbar/M^*a)$.  To see that it isolates $\tp(\cbar/M^*a)$, choose any $L(M^*)$-formula $\delta(x,\ybar)$ and suppose
$\delta(a,\cbar)$ holds.   Choose any $\dbar$ such that $\psi(a,\dbar)$ holds and it suffices to show that $\delta(a,\dbar)$ holds.
For this, we first show that $\dbar$ realizes $(p|M^*)^{(k)}$. Since $\theta(d_i)$ holds, $\tp(d_i/M^*)$ is $p$-simple of $p$-weight $\le 1$ for each $i<k$.
However, $w_p(a/M^*)=k$ and $\beta(a,\dbar)$ implies that  $w_p(a/M^*\dbar)=0$.
It follows that we must have $w_p(d_i/M^*)=1$ for each $i<k$ and moreover, $\{d_i:i<k\}$ is independent over $M^*$.  
As $p|M^*$ is strongly regular, i.e., by Theorem~\ref{olddefinability}(3), each $d_i$ realizes $p|M^*$, hence $\tp(\dbar/M^*)=(p|M^*)^{(k)}=\tp(\cbar/M^*)$.
To finish, since $\alpha(x,\cbar)$ isolates $\tp(a/M^*\cbar)$ and since $\delta(a,\cbar)$ holds, 
$\forall x(\alpha(x,\ybar)\rightarrow \delta(x,\ybar))\in \tp(\cbar/M^*)$, so  we have $\forall x(\alpha(x,\dbar)\wedge \delta(x,\dbar))$ holding.
Thus $\delta(a,\dbar)$ holds, so $\tp(\cbar/M^*a)$ is isolated.

As $\tp(a/M^*)$ is depth-zero like and DI, the same is true of $\tp(a/M)$. 
Let $N$ be any $\ell$-constructible model over $Ma$.  By Fact~\ref{lisodom}(2), $N$ is dominated by $a$ over $M$.   
So Proposition~\ref{DIprop}(4) implies that $N$ is constructible over $Ma$ and
%Thus, as any $\ell$-constructible model $N$ over $Ma$ is dominated by $a$ over $M$,  
Proposition~\ref{regmin}(2) gives the minimality of $N$ over $Ma$.
\endproof

\appendix

\section{Appendix}\label{appendix}
In this appendix, we bring together for the reader's convenience, many of the basic definitions and facts from classification theory and geometric stability theory that are used throughout the paper.  Many of these  can be found in e.g.,  Chapters 1,7, 8 of \cite{Pillay}. We assume that the reader is familiar with stability theory, independence and the basics of superstability.  Throughout this appendix, $T$ will be (at least) a stable theory.

\begin{Definition} {\em 
\begin{itemize}
\item  We say that $M$ is an {\em  $a$-model} if every strong type over every finite $B\subseteq M$  is realized in $M$. In the original notation  of Shelah (\cite{Shc}) this corresponds to ${\bf F^a_{\aleph_0}}$-saturation.
\item   Let $M\preceq N$ be models of $T$, $N$ is an {\em $na$-extension}  of $M$, denoted $M\na N$,  if for every formula $\phi(x,y)$, for every tuple $a$ from $ M$ and every finite subset $F$ of $M$, if $N$ contains a solution to $\phi(x,a)$ not in $M$, then $M$ contains a solution to $\phi(x,a)$ that  is not algebraic over $F$.
\end{itemize}}
 \end{Definition}
 
 Obviously, if $M$ is an a-model, then $M\subseteq_{na} \C$, hence $M\subseteq_{na} N$  for any $N\succeq M$.  
 We see below (e.g., Facts~\ref{qregularpair} and \ref{existdominated})
  many of the desirable attributes of working over
 $a$-models are reflected in  $na$-substructures.
The utility of this notion is that  by an easy L\"owenhiem-Skolem argument, models $M\subseteq_{na} \C$ of size
 $|T|$ exist, whereas typically
 $a$-models have larger cardinality.   
%  We see below (e.g., Facts~\ref{qregularpair} and \ref{existdominated})
%  many of the desirable attributes of working over
% $a$-models are reflected in  $na$-substructures.

\begin{Definition} {\em
\begin{itemize}
\item   If $B\subseteq A$, types  $p\in S(A)$ and $q\in S(B)$ are {\em almost orthogonal}, denoted $p\perp^a q$,  if for all  $a$ realizing $p$ and $b$ realizing $q$, 
if $\fg b B A$, then $\fg a B b$.  $p$ is almost orthogonal to the set $B$, $p\perp^a B$ if $p\perp^a q$ for every $q\in S(B)$.
\item  Two stationary types $p\in S(A)$ and $q\in S(B)$ are {\em orthogonal}, denoted $p\perp q$,  if $p|C \perp^a q$ for all $C \supseteq AB$.
$p\in S(A)$ is {\em orthogonal to a set $B$}, $p\perp B$, if $p\perp q$ for every $q\in S(B)$.

\item   A set $C$ is {\em dominated} by $E$ over $D$ if for all $b $ such that $b$ is independent from $E$ over $D$, $b$ is independent from $C$ over $D$.
\item A non-algebraic stationary type $p$ is {\em regular} if it is orthogonal to all its forking extensions.
\item  A non-algebraic stationary type $p\in S(A)$ is {\em strongly regular} if there is a formula $\phi(x)\in p$ such that, for any $B\supseteq A$, every
stationary type $q\in S(B)$ containing $\phi(x)$ is either equal to $p|B$ or orthogonal to $p$.
\end{itemize}}

\end{Definition}

Note that if $B\subseteq A$ and $p\in S(A)$, then $p\perp^a B$ if and only if $Aa$ is dominated by $A$ over $B$ for some/every $a$ realizing $p$.

We frequently use the following fact without mention (see e.g., 1.4.4.2(ii) of \cite{Pillay}).
 Because of it, many of our adjectives for regular types (e.g., local modularity, depth zero, $p$-weight)
are really features of the non-orthogonality class of the regular type.

\begin{Fact}  Non-orthogonality is an equivalence relation on set of regular types.  
\end{Fact}  

We list two existence theorems for regular types between models, which are respectively Propositions 8.3.2 and 8.3.5 of \cite{Pillay}.

\begin{Fact} \label{regularpair}  If $T$ is superstable and $M\prec N$ with $M\neq N$, then there is some $a\in N\setminus M$ with $\tp(a/M)$ regular.
\end{Fact}

\begin{Fact}  \label{qregularpair}  If $T$ is superstable,  $M\subseteq_{na} N$, and $q$ is a regular type with $\tp(N/M)\not\perp q$, then there is an $a\in N\setminus M$
with $\tp(a/M)$ regular and non-orthogonal to $q$.
\end{Fact}  

%Note that if $B\subseteq A$ and $p\in S(A)$, then $p\perp^a B$ if and only if $Aa$ is dominated by $A$ over $B$ for some/every $a$ realizing $p$.
The following fact is a consequence of Lemma 5.3 in \cite{ShelahBuechler}:

\begin{Fact} \label{existdominated} Suppose $T$ is superstable, $M\na N$ are models of $T$ and $M \subseteq A \subseteq N$. Then there exists a model $N'$, 
$M \subseteq A \subseteq N'\na N$ such that $N'$ is dominated by $A$ over $M$. 
\end{Fact}

\begin{Fact} \label{original3model}[The  3-model  lemma] {\em (Proposition 8.3.6 in \cite{Pillay})} Suppose $T$ superstable. Let $M_0\preceq M_1 \preceq M_2$ be models of $T$ such that $M_0 \na M_2$.  Suppose $a \in M_2$ and $tp(a/M_1)$ is regular non-orthogonal to $M_0$. Then there is $b \in M_2$ such that  $tp(b/M_1)$ is regular and does not fork over $M_0$ and  $b$ is not  independent  from $a$ over $M_1$. \end{Fact}

\subsection{Domination and isolation in stable theories}  \label{appendix-domisol}
In this subsection we list a number of facts about domination and isolation  that hold in arbitrary stable theories.

The first two of the following facts are obtained by forking calculus, and the third appears in 1.4.3.4 of \cite{Pillay}.

\begin{Fact}  \label{domfact}  Let $T$ be stable.
\begin{enumerate}
\item  For any  $a,b,c,D$, if $abc$ is dominated by $bc$ over $D$ and $bc$ is dominated by $c$ over $D$, then $abc$ is dominated by $c$ over $D$.
\item  If $E\supseteq D$ and $\fg {ab} D E$, then $ab$ is dominated by $b$ over $D$ if and only if $ab$ is dominated by $b$ over $E$.
%\item  If $M$ is an $a$-model and $M[b]$ is $a$-prime over $Mb$, then for any $c$, $bc$ is dominated by $b$ over $M$ if and only if  $\tp(c/Mb)$ is realized in $M[b]$.
\item  If $M$ is an $a$-model, then  $bc$ is dominated by $b$ over $M$ if and only if  $\stp(b/Mc)$ is $a$-isolated if and only if $c$ is contained in some $a$-prime model
$M[b]$ over $Mb$.
\end{enumerate}
\end{Fact}

We now recall definitions and facts about isolation and constructibility. 
  
  A type $p\in S(A)$ is {\em isolated} if there is some formula $\phi(x,a)\in p$ such that $\phi(x,a)\vdash p$.
A {\em construction sequence over $A$} is a sequence $\<a_\alpha:\alpha<\beta\>$ such that $\tp(a_\alpha/A\cup\{a_\gamma:\gamma<\alpha\})$
is isolated for every $\alpha<\beta$.  A model $N$ is {\em constructible over A} if there is a construction sequence over $A$ whose union is $N$.
If $N$ is constructible over $A$ then it is both prime and atomic over $A$. Any two constructible models over $A$ are isomorphic over $A$.
%As $T$ is countable, it follows from results of Vaught that for $A$ countable, there is a constructible model over $A$.

If $T$ is $\aleph_0$-stable, then constructible models exist over every set $A$. 
 In a superstable theory it is not always true that there are constructible models over all sets.
Indeed, one of the main goals of this paper is to determine when constructible models over  particular sets exist.

A weaker notion is {\em $\ell$-isolation}.  A type $p\in S(A)$ is {\em $\ell$-isolated} if, for every formula $\phi(x,y)$ there is a formula $\psi(x,a)\in p$
such that $\psi(x,a)\vdash p\mr{\phi}$, the restriction of $p$ to instances of $\pm\phi(x,b)$ for $b\in A$.  $\ell$-construction sequences and $N$ being
$\ell$-constructible over $A$ are defined analogously.  An advantage is that for any countable stable theory $T$,  for any set $A$ and any consistent $\phi(x,a)$
with $a\subseteq A$, an $\ell$-isolated $p\in S(A)$ extending $\phi(x,a)$ exists (see e.g., IV~2.18(4)~of~\cite{Shc}). By iterating this fact, 
$\ell$-constructible models over $A$ exist
over any  set $A$.  
The disadvantage is that there can be many non-isomorphic $\ell$-constructible models over $A$.
The following facts are well known.

%\begin{Fact}  \label{lisodom} If $T$ is stable and  $\tp(c/MA)$ is $\ell$-isolated, then $c$ is dominated by $A$ over $M$.  Iterating this gives
%that if $N$ is constructible over $MA$, then $N$ is dominated by $A$ over $M$.  
%\end{Fact}

\begin{Fact}  \label{lisodom} Suppose  $T$ is stable, $A,B$ are independent over a model $M$ and $\tp(c/MA)$ is $\ell$-isolated.
Then:
\begin{enumerate}
\item   $\tp(c/MA)\vdash\tp(c/MAB)$, so  $\tp(c/MAB)$ is $\ell$-isolated via the same formulas witnessing the $\ell$-isolation of $\tp(c/MA)$;
\item  $Ac$ is dominated by $A$ over $M$; and
\item  $\tp(c/MA)$ is isolated if and only if $\tp(c/MAB)$ is isolated.
\end{enumerate}
\end{Fact}

\bp  (1) Choose any $L$-formula $\phi(x,y)$ and assume $\theta(x,ma)\vdash\tp_\phi(c/MA)$.  If $\theta(x,ma)\not\vdash\tp_\phi(c/MAB)$, then 
$\exists x_1\exists x_2(\theta(x_1,ma)\wedge\theta(x_2,ma)\wedge\neg[\phi(x_1,mab)\leftrightarrow\phi(x_2,mab)])$.  As $\fg A M B$, finite satisfiability
would imply the existence of $b'\in M$ satisfying this, contradicting $\theta(x,ma)\vdash\tp_\phi(c/Ma)$.

(2)  Choose any $e$ with $\fg A M e$.   A non-forking extension of $\tp(c/MA)$ to $S(MAe)$ exists, and by (1) it is implied by $\tp(c/MA)$.  Thus, $\fg c {MA} e$.

(3)  Left to right is (1) and the converse is the Open Mapping Theorem.
\endproof

Iterating this yields

\begin{Fact}  \label{freeconst} 
If $T$ is stable and $A,B$ are independent over a model $M$, then for any sequence $\cbar=\<c_\alpha:\alpha<\beta\>$,
$\cbar$ is a construction sequence over $MA$ if and only if $\cbar$ is a construction sequence over $MAB$.
In particular, if $N$ is constructible over $MA$, then $N$ (as a set) is a construction sequence over $MAB$ and hence is dominated by $AB$ over $M$.
% If $M\preceq N$ are models of a stable theory, $B$ any set, and  $B\<a_\alpha:\alpha<\beta\>$ is independent from $N$ over $M$, then
%$\<a_\alpha:\alpha<\beta\>$ is a construction sequence over $MB$ if and only if it is a construction sequence over $NB$.  Similarly for $\ell$-construction sequences.
\end{Fact}

The following will be useful under the assumption of PMOP (see Section~\ref{defisolation}).

\begin{Fact}  \label{PMOPapp}  Suppose $T$ is stable, $(M_0,M_1,M_2)$ is an independent triple of models with
$b\in M_1$, and suppose that $M_1$ is constructible over $Mb$ and $M^*$ is constructible over $M_1\cup M_2$.  
If $M_1$ is dominated by $b$ over $M$, then $M^*$ is constructible over $M_2b$ and dominated by $b$ over $M_2$.  
\end{Fact}

\bp  By Fact~\ref{freeconst} the set $M_1$ is a construction sequence over $M_2b$, so since $M^*$ is constructible over $M_1M_2$, $M^*$ is constructible over $M_2b$
simply by concatenating the two construction sequences.  The domination is by Fact~\ref{lisodom}(2).
\endproof

\subsection{$p$-simplicity and locally modular regular types}\label{appendix-locallymodular}
Let $p$ be any (stationary) regular type, which for convenience we take to be over $\emptyset$.
Consider $p(\C)$ the set of realizations of the type $p$, then $(p(\C),\cl_{fork})$ forms a homogeneous 
pre-geometry (see e.g., Chapter 7  of \cite{Pillay}), where $a\in\cl_{fork}(B)$ means that $a$ forks with $B$ over $\emptyset$.
Recall that a pre-geometry $(G,\cl)$ is {\em modular} if, for all closed sets $X,Y\subseteq G$ we have
$$\dim(\cl(X\cup Y))+\dim(X\cap Y)=\dim(X)+\dim(Y)$$
and is {\em locally modular} if the above holds whenever $\dim(X\cap Y)\neq 0$.  
By considering a minimal counterexample, it is well known that  if $(G,\cl)$ is not modular, then one can find closed sets $X,Y,Z\subseteq G$, each of finite dimension, such
that $X\supseteq Z$, $Y\supseteq Z$, $\dim(X/Z)=\dim(Y/Z)=2$,  $X\cap Y=Z$, but $\dim(\cl(X\cup Y)/Z)=3$.  

%Now suppose $M$ is an a-model and $p\in S(M)$ is regular.   Then, in the pre-geometry $(p(\C),\cl_M)$ (where $a\in \cl_M(B)$ iff $a$ forks with $B$ over $M$),
%since $

Now suppose that $M$ is an a-model, let $p|M$ denote the non-forking extension of $p$ to $S(M)$.    
Since $p(M)$ contains closed  sets of infinite dimension, the induced pre-geometry $(p|M(\C),\cl_M)$, where $a\in\cl_M(B)$ iff $a$ forks with $B$ over $M$,
is locally modular if and only if it is modular.  Moreover, if 
$(p|M,\cl_M)$ is not modular, then there exist four realizations $(a_1,a_2,b_1,b_2)$ of $p|M$ such that $(a_1,a_2)$ and $(b_1,b_2)$ are each independent pairs over $M$,
$\dim(a_1a_2,b_1,b_2/M)=3$, but $\cl_M(a_1a_2)\cap\cl_M(b_1b_2)\cap p|M(\C)=\emptyset$.

%Call  a subset $X\subseteq p|M(\C)$ {\em closed} if whenever $a\in p(\C)$ does not fork over $M$,
%but $\nfg a M X$, then $a\in X$.  Because of the van der Waerden axioms, closed sets $X$ come equipped with
%a dimension $\dim(X)$, which here is equal to the cardinality of a maximal $M$-independent subset $I\subseteq M$.
%
%\begin{Definition}  {\em  The regular type $p$ is {\em locally modular} if, for any a-model $M$ and for any
%closed subsets $X,Y\subseteq (p|M)(\C)$,
%$$\dim(\cl(X\cup Y))+\dim(\cl(X\cap Y))=\dim(X)+\dim(Y)$$
%}
%\end{Definition}
%
%
%
%
%\begin{Proposition}  \label{amodelok1} \cite{Pillay}
%Suppose that $M$ is an a-model.  The following are equivalent for a regular type $p\in S(M)$:
%\begin{enumerate}
%\item  The type $p$ is not  locally modular;
%%\item  $(p(\C),\clf)$ is not modular;
%\item  There is a set of four realizations $\{a_1,a_2,b_1,b_2\}$ of $p$,
%that are dependent over $M$, yet the closures of $\{a_1,a_2\}$ and $\{b_1,b_2\}$ in $p(\C)$ are disjoint.
%\end{enumerate}
%\end{Proposition}

For many applications it is useful to work in a wider space than  $(p(\C),\clf)$. 
%Whereas the above is classical, in model theory it is useful to work in a wider space. 
We recall the definition of $p$-simplicity. 

\begin{Definition} \label{defpsimple} {\em   Suppose that a regular type $p$ is non-orthogonal to a set $A$.
\begin{itemize}
\item  A  strong type $\stp(a/A)$ is {\em hereditarily orthogonal to $p$} if  $\stp(a/B)$ is orthogonal to $p$ for every $B\supseteq A$.
\item  A strong type $\stp(a/A)$ is {\em $p$-simple} if for some a-model $M$ independent from a over $A$, 
 there is an $M$-independent set
$\{b_1,\dots,b_k\}$ of realizations of $p|M$ such that $\stp(a/Mb_1,\dots,b_k)$ is hereditarily orthogonal to $p$.  We say that  $\stp(a/A)$ is {\em $p$-simple  of weight $k$} if $k$ is least such. 
\item  If $\stp(a/A)$ is $p$-simple of $p$-weight $k$ we write $w_p(a/A)=k$.  
\item A formula $\theta(x)$ over $A$ is {\em $p$-simple of $p$-weight  $k$} if every  strong type extending $\theta$ is $p$-simple  and  $k$ is the maximum of 
$\{w_p(a/A) : \theta(x)\in\stp(a/A)\}$.
\item  For a formula $\theta(x)$ over $A$, put $(\theta_A)^{eq}:=\dcl^{eq}(A\cup\theta(\C))$.   [If $\theta(x)$ is $p$-simple, then
$\stp(b/A)$ is $p$-simple for every $b\in (\theta_A)^{eq}$.]
\item For a $p$-simple $\theta(x)$ over $A$, we say {\em $p$-weight is definable and continuous inside $\theta(x)$} 
if, for any $b\in (\theta_A)^{eq}$ and any $c\in \C^{eq}$, if $w_p(b/Ac)=n$, then there is some $\phi(x,y)\in\tp(bc/A)$ such that
for any $b'c'$ realizing $\phi$, $w_p(\phi(x,c')=n$ [hence $w_p(b'/Ac')\le n$].
%if, for all $C\supseteq A$ and for all $a$ realizing $\theta$,  if $w_p(a/C)= m$, then there is a formula $\phi(x,c)\in\tp(a/C)$ of $p$-weight
%$m$ and a formula $\psi(y)\in\tp(c/A)$ such that $w_p(\phi(x,c'))\le m$ for all $c'$ realizing $\psi$.
\end{itemize}
  
}
\end{Definition}

The following facts will be used repeatedly.  

\begin{Fact}  \label{pfacts}   Suppose $p$ is a regular type non-orthogonal to a set $A$.

\begin{enumerate}
\item  If $A\subseteq B$ and $\fg a A B$, then $\stp(a/A)$ is hereditarily orthogonal to $p$ if and only if $\stp(a/B)$ is hereditarily orthogonal to $p$.
[7.1.3' of \cite{Pillay}]
\item  In particular, (taking $B$ to be an a-model) $w_p(a/A)=0$ if and only if $\stp(a/A)$ is hereditarily orthogonal to $p$.
\item  If $\stp(a/A)$ and $\stp(b/A)$ are $p$-simple, then so is $\stp(ab/A)$.  Moreover, $w_p(ab/A)=w_p(a/Ab)+w_p(b/A)$.  [Lemmas 7.1.4(i) and 7.1.11 of \cite{Pillay}.]
\item  If $\stp(a/A)$ is $p$-simple and $b\in\acl(Aa)$, then $\stp(b/A)$ is $p$-simple.   [Lemma 7.1.4(ii) of \cite{Pillay}]
\item  If $\stp(a/A)$ is $p$-simple then $\stp(a/B)$ is $p$-simple for every $B\supseteq A$.  [Lemma 7.1.10 of \cite{Pillay}]
\item  $p$ is $p$-simple and $w_p(p)=1$.
\end{enumerate}
\end{Fact}

%As any forking extension of $p$ is hereditarily orthogonal to $p$, it follows that $p$ is $p$-simple with $w_p(p)=1$.  
%For any set $B\subseteq \C^{eq}$, put
%$\clp(B)=\{a\in \C^{\eq}:w_p(a/B)=0\}$.    Note that for any $a,B\subseteq p(\C)$, $a\in\clp(B)$ if and only if $a\in\cl(B)$ as defined above.

Non-orthogonality to $p$ gives the existence of $p$-simple types within the definable closure:

\begin{Fact} \label{psimpleindcl} (Lemma 1.17, Chapter 7 in \cite{Pillay}). Let $X$ be algebraically closed and $tp(a/X)$ be non-orthogonal to $p$. Then there is $e \in dcl(aX)$ such that $tp(e/X)$ is $p$-simple of positive $p$-weight.\end{Fact}

 Following \cite{UdiThesis} and \cite{Pillay}, for any  set $A$ with $p\not\perp A$, let 
$$D(p,A):=\{a\in\C:\stp(a/A)\ \hbox{is $p$-simple of finite $p$-weight}\}$$
We equip $D(p,A)$ with a closure operator $\clp$, namely for $a,B$ from $D(p,A)$, 
$a\in\cl_p(B)$  if and only if $w_p(a/BA)=0$.  Technically, the closure relation $\clp$ depends on $A$,
but much of the time we will take $A=\emptyset$, so we do not muddy our notation by referring to it explicitly.
In light of Fact~\ref{pfacts}, the closure space $(D(p,A),\clp)$ admits a good dimension theory, but it is not a pre-geometry as the Exchange Axiom fails.

\begin{Definition} {\em  Fix any set $A$ with $p\not\perp A$,  $D(p,A)$ is {\em modular} if $w_p(a/A)+w_p(b/A)=w_p(ab/A)+w_p((\clp(a)\cap\clp(b))/A)$ for all $a,b\in D(p,A)$.
}
\end{Definition}

As shown for example in 7.2.4 of \cite{Pillay}:

\begin{Fact}The regular type $p$ is locally modular (as defined above) if and only if  $D(p,A)$ is modular for all sets $A\not\perp p$.
\end{Fact}

 \subsection{$p$-semiregular types}\label{appendix-semiregular}
 
Within the space $D(p,A)$, it is useful to identify the $p$-semiregular types.

\begin{Definition}  {\em  $\stp(a/A)$ is {\em $p$-semiregular of weight $k$} if it is $p$-simple and is (eventually) domination equivalent to $p^{(k)}$ for some finite $k\ge 1$,
i.e., for some (equivalently, for all) a-models $M$ independent from $a$ over $A$, there is an $M$-independent sequence $\bbar=\<b_1,\dots,b_k\>$ of realizations of
the non-forking extension $p|M$ witnessing the $p$-simplicity of $\stp(a/A)$,
with $a$ and $\bbar$ domination equivalent over $M$ (for any set $X$, $\fg X M a$ if and only if $\fg X M {\bbar}$.)
}
\end{Definition}

There is a  natural Criterion for determining whether a $p$-simple type is $p$-semiregular.
This Criterion  appears as  either Fact 1.4 of \cite{HrSh} or
7.1.18 of \cite{Pillay}):

\begin{Criterion}  \label{Criterion}
 Suppose $\stp(a/X)$ is $p$-simple of positive $p$-weight, and choose
$Y\subseteq\dcl(aX)$.  Then $\stp(a/Y)$ is $p$-semiregular of positive
$p$-weight if and only if $a\not\in\acl(Y)$, but $w_p(e/Y)>0$ for every
$e\in\dcl(aY)\setminus \acl(Y)$.
\end{Criterion}

To get the existence of a $p$-semiregular type nearby a given $p$-simple type,
we  couple this with the following easy Lemma, whose proof only requires superstability.

\begin{Lemma}  \label{standard}  Suppose $a$ and $X$ are given with $a$ finite, and $Y$ is chosen arbitrarily such that $X\subseteq Y\subseteq\acl(Xa)$.
Then there is a finite sequence $b$ from $Y$ such that $Y\subseteq\acl(Xb)$.
\end{Lemma}

\bp  Recursively construct a sequence $\<b_i:i\>$ from $Y$ of maximal length
such that $\nfg a {B_i} {b_i}$, where $B_i:=X\cup\{b_j:j<i\}$.
Clearly, $R^\infty(a/B_i)$ is strictly decreasing with $i$, so any such sequence has finite length.  But, for the sequence to terminate, it must be that $Y\subseteq \acl(B_{i^*})$ for the terminal $i^*$.
\endproof

Finally, we get our existence lemma.

\begin{Lemma} \label{sr-existence}
 If $\stp(a/X)$ is $p$-simple of positive $p$-weight,
then there is a finite $b$ from $\dcl(aX)\cap \Clp(X)$ such that
$\stp(a/Xb)$ is $p$-semiregular and $w_p(a/Xb)=w_p(a/X)$.
\end{Lemma}

\bp  Let $Y=\dcl(aX)\cap\Clp(X)$ and choose a finite $b$ from $Y$ such that
$Y\subseteq\acl(Xb)$.  Now $\dcl(Ya)=\dcl(Xa)$, so if $e\in\dcl(Ya)\setminus 
\acl(Y)$, we must have $w_p(e/Y)>0$, lest we would have $e\in Y$.
Thus, Criterion~\ref{Criterion} for $p$-semiregularity applies.
\endproof

Next we record ways in which an existing $p$-semiregular type is persistent.

\begin{Lemma}  \label{semiregpersists} Suppose $\stp(e/X)$ is $p$-semiregular.
\begin{enumerate}
\item  If $e'\in\acl(eX)\setminus \acl(X)$, then $\stp(e'/X)$ is $p$-semiregular; 
\item  If $\stp(e/Y)$ is parallel to $\stp(e/X)$, then $\stp(e/Y)$ is $p$-semiregular of the same $p$-weight;
\item  If $X\subseteq Y\subseteq \Clp(X)$, then $\stp(e/Y)$ is $p$-semiregular of the same $p$-weight.
\end{enumerate}
\end{Lemma}

\bp (1) As $\dcl(e'X)\subseteq\dcl(eX)$,  the result follows by  Criterion~\ref{Criterion}.

(2)  This is immediate, as `domination equivalent to $p^{(k)}$' is preserved.

(3)  Since $\stp(e/X)$ is $p$-semiregular, we automatically have $\fg e X Y$,
so (3) follows from (2).
\endproof

Recalling Definition~\ref{defpsimple}, the  following fact is Theorem 2(b) in \cite{HrSh}. 

\begin{Fact}\label{pweightprov} Let $T$ be superstable, let $p$ be a non-trivial regular type  of depth zero and let $stp(a/B)$ be $p$-semiregular. Then $a$ lies in some $acl(B)$-definable set $D$ such that $p$-weight is continuous and definable inside $D$.\end{Fact}

\subsection{$p$-disjointness}  

\begin{Definition}  {\em  
Suppose $a,b\in D(p,C)$.  We say that {\em $a$ and $b$ are $p$-disjoint over $C$} if $\clp(Ca)\cap\clp(Cb)=\clp(C)$.
}
\end{Definition}

The next two Lemmas discuss the relationship between $p$-disjointness and forking, at least when
$\stp(ab/C)$ is $p$-semiregular.

\begin{Lemma} \label{nonforking} Suppose that  $\stp(ab/C)$ is $p$-semiregular, $C\subseteq D$ and $\stp(ab/D)$ does not fork over $C$.
Then
 $\clp(Ca)\cap\clp(Cb)=\clp(C)$ if and only if
$\clp(Da)\cap\clp(Db)=\clp(D)$.
\end{Lemma}
\bp  First, assume there is a
 `bad element' $e$ for the triple $(a,b,C)$, that is $e \in \clp(Ca)\cap\clp(Cb)\setminus \clp(C)$.  As the existence of such an $e$ is clearly determined by $\tp(ab/C)$,  by replacing
$D$ by some independent $D^*$ realizing the same strong type as $D$ over $Cab$, we may assume that $\fg {abe} C D$.
It follows immediately that 
$e\in[\clp(Da)\cap\clp(Db)]\setminus\clp(D)$ so $e$ is bad for $(a,b,D)$ as well.

Conversely, if $e$ is a `bad element' for $(a,b,D)$, let $h:=\Cb(De/Cab)$.  
We first claim that $h\not\in\clp(C)$.   If it were, then as $\stp(ab/C)$ is $p$-semiregular, we would have $\fg {ab} C h$.  But, as $\fg {ab} {Ch} {De}$, this would imply
$\fg {ab} {D} e$, contradicting $e\not\in\clp(D)$.  Thus, $h\not\in\clp(C)$. 

So, arguing by symmetry between $a$ and $b$, it suffices to prove that $h\in\clp(Ca)$.
Choose a Morley sequence $\<D_1e_1,\dots,D_ne_n\>$ in $\stp(De/Cab)$ with $D_1e_1=De$ such that $h\in dcl(e_1\ldots e_n D_1 \ldots D_n)$.
The standard argument yields
$$\fg {D_1,\dots,D_n} C {ab}$$
As well, $h\in\acl(Cab)$, hence $\fg {D_1,\dots,D_n} {Ca} h$.  Because of this, it suffices to prove that
$\stp(h/CaD_1\dots D_n)$ is hereditarily orthogonal to $p$, i.e., has $p$-weight zero.  
However, for each $i$, $w_p(e_i/aD_i)=0$, so $w_p(e_i/CaD_1\dots D_n)=0$ for each $i$.
But $h\in\dcl(e_1\dots e_n D_1\dots D_n)$, so $w_p(h/Ca D_1\dots D_n)=0$.  
\endproof

\begin{Lemma}\label{cons-p-disjoint} Suppose that $\stp(ab/C)$ is $p$-semiregular and $\clp(Ca)\cap\clp(Cb)=\clp(C)$.  
Then $\acl(Ca)\cap\acl(Cb)=\acl(C)$.
\end{Lemma}
  
\bp   Choose any $e\in\acl(Ca)\cap\acl(Cb)$.  As $\acl(Ca)\subseteq\Clp(Ca)$
and $\acl(Cb)\subseteq\Clp(Cb)$, our hypothesis implies that $e\in\Clp(C)$.
However, as $\stp(ab/C)$ is $p$-semiregular, this implies
$\fg {ab} C e$.  As $e\in\acl(abC)$, this implies $e\in\acl(C)$ as desired.
\endproof

Thanks to Lemma \ref{cons-p-disjoint} we will be able to apply the  following general result  below  about forking to $p$-disjoint $p$-semiregular
types. 

\begin{Lemma}  \label{fork}  For all $a,b,C$,  $\acl(Ca)\cap\acl(Cb)=\acl(C)$ if and only if
for every set $X$, if $\fg X {Ca} b$ and $\fg X {Cb} a$ both hold, then
$\fg X C {ab}$ holds as well.
\end{Lemma}

\bp To ease notation, assume $C=\emptyset$.  
It suffices to prove this for finite sets $X$. 
For left to right, fix an $X$, and let $D$ denote the
canonical base of $\tp(X/ab)$.  On one hand, $D\subseteq\acl(a)$, and on the
other hand, $D\subseteq\acl(b)$.  Thus, by our assumption, $D\subseteq\acl(\emptyset)$,
implying that $ X \anchor {ab}$.

For the converse,  choose any $h\in\acl(a)\cap\acl(b)$.
Then, for trivial reasons we have $\fg h a b$ and $\fg h b a$, so by our hypothesis
we have $h \anchor {ab}$.  But, as $h\in\acl(ab)$, this implies
$h\in\acl(\emptyset)$.
\endproof

\subsection{Classifiable theories}\label{defisolation}

\begin{Definition}  \label{basicdef}   {\em  Let $T$ be any stable theory.

\begin{itemize}   
\item  An {\em independent triple}  is a sequence $(A_0,A_1,A_2)$ satisfying $A_1\cap A_2=A_0$ and $\fg {A_1} {A_0} {A_2}$.  
\item   A superstable theory {\em does not have the dimensional order property} (NDOP) if, for every independent  triple $\M= (M_0,M_1,M_2)$  of a-models,
the a-prime model $M^*$ over $M_1M_2$ (which exists in any superstable theory) is minimal among all a-models containing $M_1M_2$.
\item  A superstable theory has {\em prime models over pairs} (PMOP)  if, for any independent triple $(M_0,M_1,M_2)$ of models, there is a constructible model over $M_1M_2$.  

\item A complete theory $T$ in a countable language is {\em classifiable} if $T$ is superstable, has prime models over pairs (PMOP) and does not have the dimensional order property (NDOP).
\end{itemize}
}
\end{Definition}

The following are essential facts about NDOP.   These appear as  Lemma~X~2.2 and Lemma~X~7.2  of \cite{Shc}.  The second is also Fact~8.4.5 of \cite{Pillay}.

\begin{Fact}\label{NDOP} 
\begin{itemize}  
\item $T$ has NDOP if and only if,  for $\M=(M_0,M_1,M_2)$ any independent triple of $a$-saturated models and   $M^*$  a-prime model over $M_1M_2$, any regular type $q$ non-orthogonal to $M^*$ is either non-orthogonal to $M_1$ or $M_2$.
\item  If $T$ has NDOP, then any non-trivial regular type has depth zero.
\end{itemize}
\end{Fact}

For countable, superstable theories with NDOP, the condition PMOP has many equivalents.  
As notation due to Harrington,  we say that an independent triple $(B_0,B_1,B_2)$ {\em extends}
the independent triple $(A_0,A_1,A_2)$ if $B_0\supseteq A_0$, $\fg {B_0} {A_0} {A_1A_2}$, $B_1\supseteq A_1B_0$,    $B_2\supseteq A_2B_0$, and 
$\fg {B_1} {B_0} {B_2}$.
If $(A_0,A_1,A_2)$ is any independent triple then a type $p\in S_n(A_1A_2)$ is {\em $V$-dominated} if, for every independent triple $(B_0,B_1,B_2)$ extending $(A_0,A_1,A_2)$,
and every realization $c$ of $p$, if $\fg c {A_1A_2} {B_0}$, then $\fg c {A_1A_2} {B_1B_2}$.  
The best known equivalent of PMOP is NOTOP, the negation of the Omitting Types Order Property, see e.g., Definition~XII~4.2 of \cite{Shc}.
Although it is not mentioned in the statement of Fact~\ref{PMOPequiv} below, Shelah says that $T$  has the {\em $(\aleph_0,2)$-extension property} if
there is a constructible model over any independent triple $(M_0,M_1,M_2)$ of {\em countable} models.
As the  notion of NOTOP is not used in this paper, we merely note its equivalence to the notions (2)--(4) that are used.

\begin{Fact}  \label{PMOPequiv}  The following are equivalent for a countable, superstable theory $T$ with NDOP:
\begin{enumerate}
\item  $T$ has NOTOP;
\item  $T$ has PMOP, i.e., for every independent triple $(M_0,M_1,M_2)$ of models, a constructible model $M^*$ over $M_1M_2$ exists;
\item  For every independent triple $(M_0,M_1,M_2)$ of models and every $n\ge 1$,  every $\ell$-constructible $p\in S_n(M_1M_2)$ is isolated;
\item For every independent triple $(A_0,A_1,A_2)$, every $V$-dominated $p\in S(A_1A_2)$ is isolated.
\end{enumerate}
\end{Fact}

\bp  The equivalence of (1) and (2) are Theorem XII, 4.3 and Lemma XII 6.2 of \cite{Shc}.

$(2)\Rightarrow(4)$:  This is Lemma 8.5.10 of \cite{Pillay}.

$(4)\Rightarrow(3)$:  Fix any independent triple $(M_0,M_1,M_2)$ of models and  any $\ell$-isolated $p\in S(M_1M_2)$.  To check that $p$ is $V$-dominated, it
suffices to take an independent triple  $(N_0,N_1,N_2)$ of models extending $(M_0,M_1,M_2)$.  So choose such an $(N_0,N_1,N_2)$ and choose any $c$ realizing
$p$ with $\fg c {M_1M_2} {N_0}$.     By Lemma~1.3 of \cite{Hart}, $M_1M_2\subseteq_{TV} N_1N_2$,
i.e., every $L(M_1M_2)$-formula $\phi(x)$ with a solution in $N_1N_2$ has a solution in $M_1M_2$.  Thus, by \cite{Shc}, $p\vdash \tp(c/N_1N_2)$, so in particular,
$\fg c {M_1M_2} {N_1N_2}$.  Hence $p$ is $V$-dominated, so $p$ is isolated by (4).

$(3)\Rightarrow(2)$:  First, choose  any independent triple $(M_0,M_1,M_2)$ of {\em countable} models.  
As $T$ is countable and stable, there is a countable $\ell$-constructible model $M^*$ over $M_1M_2$.  As each finite tuple from $M^*$ is $\ell$-isolated, it follows from (3) that
$M^*$ is atomic over $M_1M_2$.  As $T$ and $M^*$ are countable, this implies $M^*$ is constructible over $M_1M_2$.  Thus, in Shelah's notation, we have proved that (3) implies
that $T$ has the $(\aleph_0,2)$-extension property.  The whole of Section XII~5 of \cite{Shc}, culminating in Conclusion~5.14, is an inductive argument showing that this implies the $(\lambda, 2)$-extension property for every infinite cardinal $\lambda$, i.e., that (2) holds.  An alternate approach to this inductive argument is given in Theorem~3.5 of \cite{Hart}.
\endproof

%As $p$ is $\ell$-isolated, there is an $\ell$-constructible model $M^*$ over $M_1M_2$ containing $c$.  By Lemma~\ref{1.3} of \cite{Hart}, $M_1M_2\subseteq_{TV} N_1N_2$,
%i.e., every $L(M_1M_2)$-formula $\phi(x)$ with a solution in $N_1N_2$ has a solution in $M_1M_2$.  Thus, by \cite{Shc}, $p\vdash \tp(c/N_1N_2)$.

\begin{Remark}  {\em  In fact, under the hypotheses of Fact~\ref{PMOPequiv}, if follows from Theorem XII 4.17 of \cite{Shc} that any constructible model $M^*$ occurring in the definition of PMOP is also minimal over
$M_1M_2$.
}
\end{Remark}

%In addition to these three equivalents that are used in this paper, the best known equivalent is with NOTOP, see e.g., Definition~\ref{ZZZ} of \cite{Shc}.
%As the notion is not used here, we simply state the following fact.
%
%\begin{Fact}  The following are equivalent for any 

\end{document}